\documentclass{amsart}

\RequirePackage{float}% figure placement
\usepackage{amsmath,amsthm,amssymb,amsfonts}
\usepackage[colorlinks=true]{hyperref}
\usepackage{graphicx}
\usepackage{tikz}
\usetikzlibrary{arrows.meta}
\usepackage{pgfplots}
\pgfplotsset{compat=1.18}
\hypersetup{urlcolor=blue, citecolor=blue, linkcolor=blue}

\newcommand{\de}{\delta}
\newcommand{\ep}{\varepsilon}
\newcommand{\vep}{\varepsilon}

\newcommand{\R}{\mathbb{R}}
\DeclareMathOperator{\supp}{supp}
\newcommand{\les}{\lesssim}

\newcommand{\la}{\lambda}
\newcommand{\beeq}{\begin{equation}}
\newcommand{\eneq}{\end{equation}}

\numberwithin{equation}{section}
\newtheorem{theorem}{Theorem}[section]
\newtheorem{lemma}[theorem]{Lemma}

\newtheorem{remark}[theorem]{Remark}

\title[CRITICAL CURVE FOR MULTI-SPEED WAVE SYSTEMS]
{Critical Curve for Three-Dimensional Radial Wave Systems with Distinct Speeds}

\dedicatory{Dedicated to Professor Hideo Kubo on the occasion of his 60th birthday.}

\author{Wei Dai}
\address{School of Mathematical Sciences\\                	
Zhejiang University of Technology\\                Hangzhou 310023, P. R. China}
\email{daiw23@zjut.edu.cn}
\author[Mengyun Liu]{Mengyun Liu\textsuperscript{*}}
\thanks{\textsuperscript{*}Corresponding author.}
\address{Department of Mathematics\\                	
Zhejiang Sci-Tech University\\                Hangzhou 310018, P. R. China}
\email{mengyunliu@zstu.edu.cn}
\author{Chengbo Wang}
\address{School of Mathematical Sciences\\
 Zhejiang Normal University\\ 
  Jinhua 321004, P. R. China}
 \email{wangcbo@zjnu.edu.cn}

 \thanks{The first author was supported by NSFC 12501292. The second author was supported by NSFC 12101558. The third author was supported by NSFC 12671293.}

\begin{document}
\begin{abstract}
We study radially symmetric solutions to the three-dimensional semilinear wave system 
with distinct propagation speeds
\[
  u_{tt}-c_1^2\Delta u=|v_t|^p\ , \ 
  v_{tt}-c_2^2\Delta v=|u_t|^q\ ,
  \   c_1, c_2>0\ , c_1\neq c_2\ .
\]
For \(1<p,q<2\), we prove that
\[
  \Gamma(p,q)
  =
  2+\max\left\{\frac1p,\frac1q\right\}-pq=0
\]
is the exact critical curve separating small-data global existence
from finite-time blow-up. More precisely, sufficiently small,
compactly supported radial initial data give rise to global solutions
when \(\Gamma(p,q)<0\). On the other hand, for a class of initial data, the corresponding
solutions blow up in finite time, with an upper bound for the lifespan,
when \(\Gamma(p,q)\geq 0\).
\end{abstract}

\maketitle

\tableofcontents

\section{Introduction}

We consider the following three-dimensional semilinear wave system with
distinct propagation speeds:
\begin{equation}
\label{eq:main-system}
\left\{
\begin{aligned}
&u_{tt}-c_1^2\Delta u
  =\lvert v_t\rvert^p,
  &&(t,x)\in[0,T)\times\mathbb R^3,\\
&v_{tt}-c_2^2\Delta v
  =\lvert u_t\rvert^q,
  &&(t,x)\in[0,T)\times\mathbb R^3,\\
&u(0,x)
  =\ep f_1(x),
  &&u_t(0,x)=\ep g_1(x),\\
&v(0,x)
  =\ep f_2(x),
  &&v_t(0,x)=\ep g_2(x).
\end{aligned}
\right.
\end{equation}
Here
\(c_1,c_2>0\) and \(\ep>0\) is a small parameter measuring the size of the initial data.
We assume that
\beeq\label{initial-data}
  (f_i,g_i)\in C_c^2(\mathbb R^3)\times C_c^1(\mathbb R^3),
  \qquad i=1,2\ .
\eneq
%and, for some fixed \(R>0\),
%\[
 % \supp f_i\cup\supp g_i
 % \subset B_R:=\{x\in\mathbb R^3:|x|\leq R\},
 % \qquad i=1,2.
%\]

One of the motivations for studying this problem comes from the work of
Kubo, Kubota, and Sunagawa \cite{KKS}. They proved the global existence
of small-data radial solutions in three space dimensions with distinct
propagation speeds under the condition
\[
  \max\left\{
    \frac{q}{pq-1},
    \frac{p}{pq-1}
  \right\}
  -1<0.
\]
By contrast, when the propagation speeds coincide, the critical curve
in three space dimensions is given by
\[
  \max\left\{
    \frac{q+1}{pq-1},
    \frac{p+1}{pq-1}
  \right\}
  =1;
\]
see \cite{KKS, Deng99, I-S-W}. It is therefore natural to ask whether one can identify the exact critical
curve for the system with distinct propagation speeds. This is precisely
the question addressed in the present paper.

Before stating our main result, we briefly review some related results on
the scalar Glassey problem and weakly coupled wave systems.

For comparison, we recall the classical scalar Glassey problem, which has
been extensively studied. For the scalar equation in \(n\) space dimensions, the global existence
of small-data solutions for
\[
  p>p_G(n):=1+\frac{2}{n-1}\ ,
\]
is known in dimensions \(n=2,3\), through the works of Sideris
\cite{sideris}, Hidano--Tsutaya \cite{hk}, and Tzvetkov \cite{nt}. See also John \cite{FJ1981} and John--Klainerman
\cite{JK84}.
On the other hand, for \(1<p\leq p_G(n)\) and \(n\geq2\), finite-time
blow-up is known, together with an upper bound for the lifespan, at least when the initial velocity is nontrivial; see Zhou \cite{zhou}
and the references therein. In higher space dimensions, the corresponding global existence problem
remains open in general, except in the radial setting; see
Hidano--Wang--Yokoyama \cite{hyw}.

Turning from the scalar equation to systems with equal propagation
speeds, Deng \cite{Deng99} obtained finite-time blow-up results in three space dimensions under the condition
\[
  q(p-1)\leq2\ , \qquad 1<p\leq q\ .
\]
On the global-existence side, Kubo, Kubota, and Sunagawa \cite{KKS}
proved the global existence of radial solutions for
sufficiently small initial data in three space dimensions when
\[
  q(p-1)>2,
  \qquad 1<p\leq q\ .
\]
Kubo and Takaki \cite{KT2010} subsequently removed the radial symmetry
assumption in the special case \(p=2\), \(q>2\), and \(n=3\).
Consequently, in the three-dimensional radial setting, the results of
Deng \cite{Deng99} and Kubo, Kubota, and Sunagawa \cite{KKS} identify
\[
  q(p-1)=2,
  \qquad 1<p\leq q,
\]
as the critical curve for the equal-speed system. For general space dimensions \(n\geq2\), Ikeda, Sobajima, and Wakasa
\cite{I-S-W} established finite-time blow-up results in the region
\[
  \Gamma_1(p,q,n)
  :=
  \max\left\{
    \frac{p+1}{pq-1},
    \frac{q+1}{pq-1}
  \right\}
  -\frac{n-1}{2}
  \geq0,
  \qquad
  p,q>1.
\]
For \(n=3\) and \(1<p\leq q\), this condition is equivalent to
$$
  q(p-1)\leq2\ .
$$

As discussed above, the situation changes substantially when
\(c_1\neq c_2\). In three space dimensions, Kubo, Kubota, and
Sunagawa \cite{KKS} established the global existence of radially
symmetric small-data solutions for the system with distinct propagation
speeds under the condition
\begin{equation*}
  \Gamma_2(p, q):=\max\left\{
    \frac{q}{pq-1},
    \frac{p}{pq-1}
  \right\}
  -1
  <0.
\end{equation*}
On the other hand, much less is known about finite-time blow-up for
systems with distinct propagation speeds. Xu \cite{X} proved finite-time
blow-up under the condition
\begin{equation*}
  \frac{p+q+2}{pq-1}\geq2n.
\end{equation*}
Thus, there remains a substantial gap between the known global-existence
and blow-up regions. 

The purpose of the present paper is to close this gap in the
three-dimensional radial setting. A central feature of
our analysis is the emergence of the curve
\[
  \Gamma(p,q)
  :=
  2+\max\left\{\frac1p,\frac1q\right\}-pq.
\]
\begin{figure}[H]
  \centering
  \includegraphics[width=0.6\textwidth]{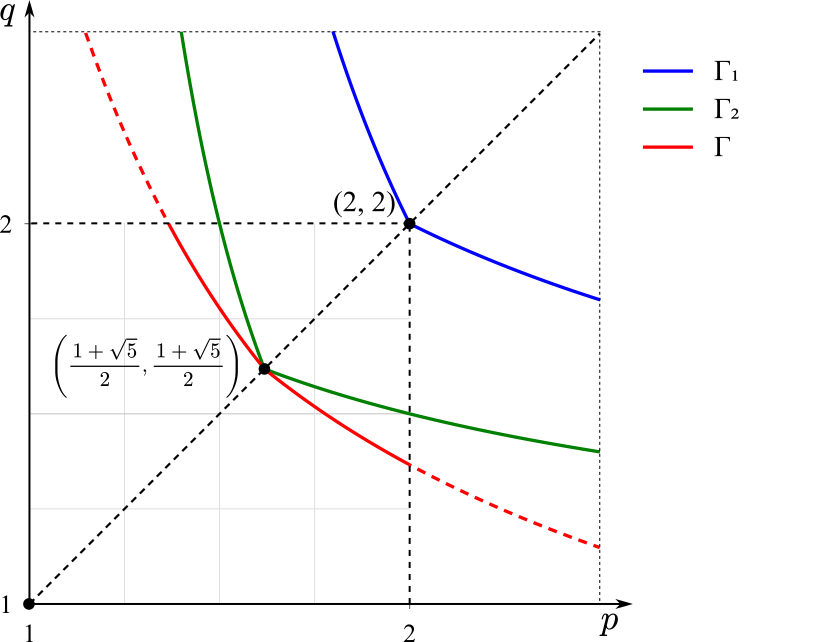}
\caption{}
  \label{fig:critical-curves}
\end{figure}
We show that, in the range \(1<p,q<2\), the zero set
\[
  \Gamma(p,q)=0
\]
is precisely the critical curve separating small-data global existence
from finite-time blow-up for compactly supported radial data. 
Figure~\ref{fig:critical-curves} compares this curve with the previously
known thresholds in the \(p\)-\(q\) plane in three space dimensions.
Here \(\Gamma_1=0\) denotes the critical curve for the corresponding
system with equal propagation speeds, while \(\Gamma_2=0\) is the
threshold obtained in the known global-existence result for systems
with distinct propagation speeds. 
In the range \(1<p,q<2\), the new curve \(\Gamma=0\) lies below
\(\Gamma_1=0\) and below \(\Gamma_2=0\), except at
\(p=q=(1+\sqrt5)/2\), where it meets \(\Gamma_2=0\).

We are now in a position to state our main result.
\subsection{Main result}

\begin{theorem}\label{thm:main}
Consider the system \eqref{eq:main-system} with the initial data
\eqref{initial-data} and $c_1\neq c_2$. Then the following assertions hold.

\begin{enumerate}
\item \emph{Global existence.}
Assume that
\[
  \Gamma(p,q)<0\ ,
  \qquad
  p,q\in(1,2).
\]
Suppose that the initial data are radially symmetric.
Then there exists \(\vep_0>0\) such that, for
\(0<\vep\leq\vep_0\), the problem \eqref{eq:main-system} admits a
global radial \(C^1\)-solution.

\item \emph{Blow-up.}
Assume that
\[
  \Gamma(p,q)\geq0\ ,
  \qquad
  p,q>1.
\]
Suppose that the initial data are radially symmetric and
\[
 f_i=0,\qquad g_i\geq0,\qquad g_i\not\equiv0,
 \qquad i=1,2.
\]
Then the corresponding radial solution
cannot exist globally. More precisely, there exist constants
\(C>0\) and \(\vep_0>0\) such that, for \(0<\vep\leq\vep_0\),
\begin{equation}\label{eq:compact-blowup-lifespan}
  T_\vep\leq
  \begin{cases}
    C\vep^{-\frac{pq-1}{\Gamma(p,q)}},
      & \Gamma(p,q)>0,\\[2mm]
    \exp\left(C\vep^{-pq\min\{p,q\}}\right),
      & \Gamma(p,q)=0\ , 
  \end{cases}
\end{equation}
where \(T_\vep\) denotes the lifespan of the corresponding radial solution.
\end{enumerate}
\end{theorem}
\begin{remark}
The compact-support assumption in the global-existence part of the theorem
is not essential. In fact, the same argument applies to sufficiently small
radial initial data satisfying the scaling-invariant decay rates.  See  Section \ref{section-5} for the
precise assumptions and the corresponding global-existence result.
\end{remark}

The main point in the global-existence argument is a sharper use of the
separation between the two characteristic cones. The estimates in
\cite{KKS} effectively rely on the integrability of one of the
characteristic weights. In the range considered here, however, the two
relevant weights may both be nonintegrable. By working directly with
the one-dimensional representation formula, we exploit the separation
of their characteristic centers and show that their combined decay is
still integrable. The resulting convolution estimate in
Section~\ref{Convolution}, together with the scaling-critical weights, closes the fixed-point argument precisely
in the range \(\Gamma(p,q)<0\).

To prove finite-time blow-up, we use the standard spatial-functional method. 
A key additional ingredient is to exploit the fundamental solution of the
three-dimensional radial wave equation to improve the lower bounds for the
spatial functionals. We then apply a Kato-type lemma to obtain finite-time blow-up and the corresponding upper bounds for the lifespan. On the critical
curve \(\Gamma(p,q)=0\), the power-type lower bounds alone are no
longer sufficient. The main new ingredient is to obtain an additional
logarithmic gain by exploiting the geometry of the two distinct
characteristic cones. Combined with a critical Kato-type argument,
this yields finite-time blow-up in the critical case.

To prepare for the proof, we first introduce the radial reduction and
the integral formulation that will be used throughout the paper.

\subsection{Radial reduction and integral formulation}

We restrict ourselves to radially symmetric solutions,
\[
  u(t,x)=u(t,\lvert x\rvert)=u(t,r),
  \qquad
  v(t,x)=v(t,\lvert x\rvert)=v(t,r),
  \qquad r=\lvert x\rvert .
\]
Set
\[
  U(t,r)=ru(t,r),
  \qquad
  V(t,r)=rv(t,r).
\]
Then the radial system \eqref{eq:main-system} is transformed into the
one-dimensional system on the half-line
\begin{equation}\label{eq:radial-UV-system}
  \begin{cases}
    U_{tt}-c_1^2 U_{rr}
    =
    r^{1-p}\lvert V_t\rvert^p,
    \\[1mm]
    V_{tt}-c_2^2 V_{rr}
    =
    r^{1-q}\lvert U_t\rvert^q.
  \end{cases}
\end{equation}
We extend radial functions evenly in the radial variable. Thus
\[
  u(t,-r)=u(t,r),
  \qquad
  v(t,-r)=v(t,r),
  \qquad r>0.
\]
Consequently, \(U=ru\) and \(V=rv\) are extended as odd functions of \(r\). 
Throughout the paper, whenever
\(r^{1-p}|V_t|^p\) and \(r^{1-q}|U_t|^q\) are used on the whole line,
they denote the odd extensions of the corresponding functions defined
for \(r>0\).
The corresponding initial data are
\[
  U(0,r)=\vep r f_1(r),
  \qquad
  U_t(0,r)=\vep r g_1(r),
\]
and
\[
  V(0,r)=\vep r f_2(r),
  \qquad
  V_t(0,r)=\vep r g_2(r).
\]
Then, after the odd extension in \(r\),
\begin{equation*}
  r f_i(r)\in C^2_{c}(\mathbb R),
  \qquad
  r g_i(r)\in C^1_{c}(\mathbb R),
  \qquad i=1,2.
\end{equation*}
Define the one-dimensional Duhamel operator
\beeq\label{Duhamel operator}
   L_{i}(F)(t,r)
  =
  \frac1{2c_i}
  \int_0^t
  \int_{r-c_i(t-s)}^{r+c_i(t-s)}
  F(s,\rho)\,d\rho ds\ .
\eneq
Then the integral formulation of \eqref{eq:radial-UV-system} is
\[
  U
  =
  U^0+L_1\bigl(r^{1-p}\lvert V_t\rvert^p\bigr),
  \qquad
  V
  =
  V^0+L_2\bigl(r^{1-q}\lvert U_t\rvert^q\bigr),
\]
where \(U^0\) and \(V^0\) denote the corresponding free one-dimensional wave
solutions. 
By a radial \(C^1\)-solution on \([0,T)\), we mean a radially symmetric pair
\[
  (u,v)\in C^1\bigl([0,T)\times\mathbb R^3\bigr)
\]
satisfying the above integral formulation. We denote by \(T_\vep\) the
maximal existence time of the corresponding radial \(C^1\)-solution.

Throughout the paper, \(A\simeq B\) means that there exist positive
constants \(c\) and \(C\), independent of \(\vep\), such that
\[
  cB\leq A\leq CB .
\]
We write \(A\lesssim B\) (respectively, \(A\gtrsim B\)) for a one-sided
upper (respectively, lower) bound up to such a constant.

The paper is organized as follows. In Section~\ref{Convolution}, we establish
the one-dimensional convolution estimates that capture the interaction
between the two characteristic weights. In Section~\ref{section-5}, we prove global existence for sufficiently small radial initial data satisfying the scaling-critical decay rates; see Theorem~\ref{thm:global-critical-decay}. In Section~\ref{sec:finite-time-blow-up}, we prove
finite-time blow-up. We first recall the subcritical Kato-type lemma and establish a critical version, then derive the spatial integral estimates needed in
the subcritical case, and finally treat the critical case by obtaining
the logarithmic improvement described above. In Appendix \ref{app:scaling-decay}, we revisit the scaling of the system and show that the
decay rates used in the global-existence argument are sharp: for either initial velocity component, any positive loss in its scaling-critical decay rate can lead to finite-time blow-up for suitable initial data. 
\section{Convolution Estimates}\label{Convolution}

The global-existence argument relies on a precise estimate for the
interaction of the two characteristic weights associated with the
distinct propagation speeds. After the radial reduction, the Duhamel
terms naturally lead to one-dimensional integrals involving these two
weights. The main point is the case
\[
  0<\beta,\gamma<1,
  \qquad
  \beta+\gamma>1,
\]
where neither characteristic weight is integrable by itself. When
\(a\ne1\), however, the separation of the corresponding characteristic
directions yields sufficient combined decay to make their product
integrable. The following convolution estimate quantifies this effect
and will be the main tool used in Section~\ref{section-5}.

\begin{lemma}[One-dimensional convolution estimate]
\label{lem:oned-convolution}
Assume that
\[
  \beta>0, \qquad \gamma>0,\qquad
  a>0,\qquad
  a\ne1.
\]
For $(t,r)\in(0,\infty)\times\mathbb R$, 
define
\[
  I(\beta,\gamma,a)
  =
  \int_0^t
  (1+|r-t+s|)^{-\beta}
  (1+||r-t+s|-as|)^{-\gamma}
  \,ds.
\]
Then
\[
  I(\beta,\gamma,a)
  \lesssim
  \begin{cases}
    (1+t+|r|)^{1-\beta-\gamma},
    & \beta<1,\ \gamma<1,\ \beta+\gamma<1,\\[0.5em]
    \log(2+t+|r|),
    & \beta<1,\ \gamma<1,\ \beta+\gamma=1,\\[0.5em]
    (1+\lvert t-r\rvert)^{1-\beta-\gamma},
    & \beta<1,\ \gamma<1,\ \beta+\gamma>1,\\[0.5em]
    (1+\lvert t-r\rvert)^{-\beta},
    & \beta\leq1,\ \gamma>1,\\[0.5em]
    (1+\lvert t-r\rvert)^{-\gamma},
    & \beta>1,\ \gamma\leq1,\\[0.5em]
    (1+\lvert t-r\rvert)^{-\min\{\beta,\gamma\}},
    & \beta>1,\ \gamma>1,\\[0.5em]
    (1+\lvert t-r\rvert)^{-\min\{\beta,\gamma\}}
    \log(2+\lvert t-r\rvert),
    & \max\{\beta,\gamma\}=1,\quad \beta+\gamma>1 .
  \end{cases}
\]
In the borderline case \(0<\beta,\gamma<1\) and
\(\beta+\gamma=1\), one has the refined bound
\begin{equation}
\label{eq:refined-critical-convolution}
  I(\beta,\gamma,a)
  \lesssim
  1+
  \log\left(
    \frac{2+t+|r|}{1+|t-r|}
  \right).
\end{equation}
\end{lemma}

\begin{remark}
For \(\gamma>1\), the estimates above follow from
Lemma~3.1 of Kubo, Kubota, and Sunagawa \cite{KKS}.
The main new case here is
\[
  0<\beta,\gamma<1,\qquad \beta+\gamma>1,
\]
where neither characteristic weight is integrable by itself.
\end{remark}

\begin{proof}
We give the details for \(a>1\). The case \(0<a<1\) follows from the
same decomposition and is discussed at the end of the proof. Set
\[
  H_\sigma(A)
  :=
  \begin{cases}
    (1+A)^{1-\sigma}, & 0<\sigma<1,\\[0.4em]
    \log(2+A), & \sigma=1,\\[0.4em]
    1, & \sigma>1.
  \end{cases}
\]

We first consider the case \(t>r\). Set
\[
  A=t-r,\qquad
  s_1=\frac{A}{a+1},\qquad
  s_2=A,\qquad
  d=s_2-s_1=\frac{a}{a+1}A.
\]
We decompose \([0,t]\) into
\[
  \begin{aligned}
  E_1
  &=
  \left\{
    s\in[0,t]:
    |s-s_1|\leq\frac d4
  \right\},  \qquad  E_2
  =
  \left\{
    s\in[0,t]:
    |s-s_2|\leq\frac d4
  \right\},\\
  E_3
  &=
  \left\{
    s\in[0,t]\cap(s_1,s_2):
    |s-s_1|\geq\frac d4,\
    |s-s_2|\geq\frac d4
  \right\}, 
  \end{aligned}
\]
and \(E_4
  =
  [0,t]\setminus(E_1\cup E_2\cup E_3).\)
  Writing \(I_j\) for the integral over \(E_j\), we have
\[
  I(\beta,\gamma,a)
  \leq I_1+I_2+I_3+I_4.
\]
For \(0\leq s\leq A\), we have
\begin{align*}
  |r-t+s|=A-s\ , 
  \bigl||r-t+s|-as\bigr|
  =
  |A-(a+1)s|
  =
  (a+1)|s-s_1|.
\end{align*}
For \(s\geq A\), we have
\begin{align*}
  |r-t+s|=s-A\ , \quad
  \bigl||r-t+s|-as\bigr|
  =
  A+(a-1)s.
\end{align*}
On \(E_1\), we have \(s<A\) and
\[
  1+|s-A|\simeq1+A.
\]
Then
\[
  \begin{aligned}
  I_1
  &\lesssim
  (1+A)^{-\beta}
  \int_{|s-s_1|\leq d/4}
  (1+|s-s_1|)^{-\gamma}\,ds\\
  &\lesssim
  (1+A)^{-\beta}
  \int_0^A(1+y)^{-\gamma}\,dy
  \lesssim
  (1+A)^{-\beta}H_\gamma(A).
  \end{aligned}
\]
We next consider \(E_2\). If \(s\in E_2\cap[0,A]\), then
\[
  \bigl||r-t+s|-as\bigr|
  =
  |A-(a+1)s|
  \simeq A.
\]
If \(s\in E_2\cap[A,t]\), then
\[
  \bigl||r-t+s|-as\bigr|
  =
  A+(a-1)s
  \simeq A.
\]
Thus, throughout \(E_2\),
\[
  1+\bigl||r-t+s|-as\bigr|
  \simeq1+A,
  \qquad
  |r-t+s|=|s-A|.
\]
It follows that
\[
  \begin{aligned}
  I_2
  &\lesssim
  (1+A)^{-\gamma}
  \int_{|s-A|\leq d/4}
  (1+|s-A|)^{-\beta}\,ds\\
  &\lesssim
  (1+A)^{-\gamma}
  \int_0^A(1+y)^{-\beta}\,dy
  \lesssim
  (1+A)^{-\gamma}H_\beta(A).
  \end{aligned}
\]
On \(E_3\), the distance from \(s\) to each of \(s_1\) and \(s_2\) is
comparable to \(A\). Hence
\[
  1+|s-A|\simeq1+A\ , \quad
  1+|A-(a+1)s|\simeq1+A.
\]
Since \(|E_3|\lesssim A\), we obtain
\[
  I_3
  \lesssim
  (1+A)^{1-\beta-\gamma}.
\]

It remains to estimate \(I_4\). The part to the left of \(E_1\) is
nonempty only when \(s_1-d/4>0\). On this interval,
\[
  1+A-s\simeq1+A\ , \quad
  1+|A-(a+1)s|\simeq1+A.
\]
Its contribution is therefore bounded by
\[
  C(1+A)^{1-\beta-\gamma}.
\]
For the part of \(E_4\) to the right of \(E_2\), we have
\[
  s>A+\frac d4.
\]
Setting \(u=s-A\), and recalling that \(d\simeq A\), we obtain
\[
  u>\frac d4\simeq A.
\]
Consequently,
\[
  1+A+(a-1)s
  =
  1+aA+(a-1)u
  \simeq
  1+u.
\]
Therefore
\[
  \begin{aligned}
  &\int_{\{s>A+d/4\}}
  (1+s-A)^{-\beta}
  \bigl(1+A+(a-1)s\bigr)^{-\gamma}\,ds\\
  &\qquad\lesssim
  \int_{d/4}^{t-A}
  (1+u)^{-\beta-\gamma}\,du\\
  &\qquad\lesssim
  \begin{cases}
    (1+t+|r|)^{1-\beta-\gamma},
    & \beta+\gamma<1,\\[0.4em]
    \log(2+t+|r|),
    & \beta+\gamma=1,\\[0.4em]
    (1+A)^{1-\beta-\gamma},
    & \beta+\gamma>1.
  \end{cases}
  \end{aligned}
\]
Here and below, an integral over an empty interval is understood to be
zero. Together with the estimates for \(I_1\), \(I_2\), and \(I_3\), this gives
the bounds in the statement. In particular, when \(\beta=1\) and
\(\gamma>1\),
\[
  (1+A)^{-\gamma}H_1(A)
  \lesssim
  (1+A)^{-1},
\]
while, when \(\beta>1\) and \(\gamma=1\),
\[
  (1+A)^{-\beta}H_1(A)
  \lesssim
  (1+A)^{-1}.
\]
Thus the result holds when \(t>r\).

We now consider the case \(r>t\). Set
\[
  A=r-t.
\]
Then
\[
  |r-t+s|=A+s,
  \qquad
  \bigl||r-t+s|-as\bigr|
  =
  |A-(a-1)s|.
\]
Let
\[
  s_0=\frac{A}{a-1}
\]
and decompose
\[
  F_1
  =
  \left[\frac{s_0}{2},2s_0\right]\cap[0,t],
  \qquad
  F_2=[0,t]\setminus F_1.
\]
On \(F_1\), we have \(s\simeq A\), and hence
\[
  1+A+s\simeq1+A.
\]
Therefore
\[
  \begin{aligned}
  &\int_{F_1}
  (1+A+s)^{-\beta}
  (1+|A-(a-1)s|)^{-\gamma}\,ds\\
  \lesssim&
  (1+A)^{-\beta}
  \int_0^{CA}(1+y)^{-\gamma}\,dy
  \lesssim
  (1+A)^{-\beta}H_\gamma(A).
  \end{aligned}
\]
On \(F_2\), either \(s<s_0/2\) or \(s>2s_0\). If \(s<s_0/2\), then
\[
  A-(a-1)s\geq\frac A2
  \simeq A+s.
\]
If \(s>2s_0\), then
\[
  (a-1)s-A\geq\frac{a-1}{2}s
  \simeq A+s.
\]
Consequently,
\[
  \begin{aligned}
  \int_{F_2}
  &(1+A+s)^{-\beta}
  (1+|A-(a-1)s|)^{-\gamma}\,ds\\
  &\lesssim
  \int_0^t(1+A+s)^{-\beta-\gamma}\,ds\\
  &\lesssim
  \begin{cases}
    (1+t+|r|)^{1-\beta-\gamma},
    & \beta+\gamma<1,\\[0.4em]
    \log(2+t+|r|),
    & \beta+\gamma=1,\\[0.4em]
    (1+A)^{1-\beta-\gamma},
    & \beta+\gamma>1.
  \end{cases}
  \end{aligned}
\]
Combining this with the estimate on \(F_1\), we obtain the bounds in the
statement. In particular, the logarithmic factor in \(H_1(A)\) is absorbed
by the stronger algebraic decay when either \(\beta=1<\gamma\) or
\(\gamma=1<\beta\).

Finally, suppose that \(r=t\). Then
\[
  |r-t+s|=s,
  \qquad
  \bigl||r-t+s|-as\bigr|
  =
  |1-a|s.
\]
Since \(a\ne1\), it follows that
\[
  \begin{aligned}
  I(\beta,\gamma,a)
  &\lesssim
  \int_0^t(1+s)^{-\beta-\gamma}\,ds
 \lesssim
  H_{\beta+\gamma}(t).
  \end{aligned}
\]
Equivalently,
\[
  I(\beta,\gamma,a)
  \lesssim
  \begin{cases}
    (1+t+|r|)^{1-\beta-\gamma},
    & \beta+\gamma<1,\\[0.4em]
    \log(2+t+|r|),
    & \beta+\gamma=1,\\[0.4em]
    1,
    & \beta+\gamma>1.
  \end{cases}
\]
Since \(|t-r|=0\), these are precisely the corresponding bounds in the
statement.

It remains to prove the refined estimate. Assume that
\(0<\beta,\gamma<1\) and \(\beta+\gamma=1\). If \(t>r\), the preceding
estimates show that \(I_1+I_2+I_3\), together with the left-hand
contribution from \(E_4\), is bounded uniformly.
For the right-hand part of \(E_4\), we retain the lower limit obtained
above and find
\[
  \begin{aligned}
  \int_{d/4}^{t-A}\frac{du}{1+u}
  &\lesssim
  1+
  \log\left(
    \frac{2+t+|r|}{1+A}
  \right).
  \end{aligned}
\]
Here we used \(d\simeq A\), and the integral is understood to be zero
when its interval is empty. Since \(A=t-r=|t-r|\), this proves
\eqref{eq:refined-critical-convolution} in this case.

If \(r>t\), the contribution of \(F_1\) is bounded uniformly, while
\[
  \begin{aligned}
  \int_{F_2}
  (1+A+s)^{-\beta}
  (1+|A-(a-1)s|)^{-\gamma}\,ds
  &\lesssim
  \int_0^t\frac{ds}{1+A+s}
  \\
  &\lesssim
  1+
  \log\left(
    \frac{2+t+|r|}{1+A}
  \right).
  \end{aligned}
\]
Now \(A=r-t=|t-r|\). Finally, when \(r=t\), the bound follows from
\(H_1(t)=\log(2+t)\). This proves
\eqref{eq:refined-critical-convolution} for \(a>1\).

It remains to consider \(0<a<1\). The argument is the same, with the
locations of the characteristic centers interchanged. Indeed, when
\(t>r\) and \(A=t-r\), the possible centers are
\[
  s_1=\frac{A}{1+a},
  \qquad
  s_2=\frac{A}{1-a},
\]
whenever they belong to \([0,t]\). When \(r>t\), one has
\[
  \bigl||r-t+s|-as\bigr|
  =
  A+(1-a)s,
\]
so no additional center occurs. Decomposing the integration interval
into neighborhoods of the possible centers and their complement gives
the same estimates as above, including
\eqref{eq:refined-critical-convolution}. Thus the conclusion also holds for
\(0<a<1\).

\end{proof}

\section{Global Existence}\label{section-5}

This section proves the global-existence part of
Theorem~\ref{thm:main}. We consider the radially symmetric reduction of
\eqref{eq:main-system}, together with its odd extension to the whole line,
namely system \eqref{eq:radial-UV-system} for \((U,V)\).

\begin{theorem}[Global existence for scaling-critical data]
\label{thm:global-critical-decay}
Assume that
\beeq\label{global-index}
 \Gamma(p,q)<0\ , \quad 1<p,q<2\ . 
\eneq
Let the radial initial profiles satisfy
\(
 (f_i,g_i)\in C^2(\mathbb R^3)\times C^1(\mathbb R^3), i=1,2,
\)
and
\begin{equation}
\left\{
\begin{aligned}
 \left\|
 (1+r)^{\alpha_*}
 \bigl(
 |(rf_1)'(r)|
 +|(rf_1)''(r)|
 +|rg_1(r)|
 +|(rg_1)'(r)|
 \bigr)
 \right\|_{L^\infty_r}
 &<\infty,
 \\
 \left\|
 (1+r)^{\beta_*}
 \bigl(
 |(rf_2)'(r)|
 +|(rf_2)''(r)|
 +|rg_2(r)|
 +|(rg_2)'(r)|
 \bigr)
 \right\|_{L^\infty_r}
 &<\infty,
\end{aligned}
\right.
\label{ini}
\end{equation}
where
\begin{equation}
\label{eq:scaling-indices}
 \alpha_*(p,q)=\frac{p+2-pq}{pq-1},
 \qquad
 \beta_*(p,q)=\frac{q+2-pq}{pq-1}.
\end{equation}
Then there exists \(\vep_0>0\) such that, for
\(0<\vep\leq\vep_0\), the problem \eqref{eq:main-system}
admits a global radial \(C^1\)-solution.
\end{theorem}

\begin{remark}
The exponents \(\alpha_*\) and \(\beta_*\) are dictated by the natural
scaling of \eqref{eq:main-system} and give the scaling-critical decay rates
of the initial data. Moreover, both rates are sharp: for either initial velocity component, any positive loss in its scaling-critical decay rate can lead to finite-time blow-up for suitable initial data. See Appendix~\ref{app:scaling-decay} for the scaling
argument and the sharpness result.
\end{remark}

\begin{proof}
After rescaling space and relabeling the two components if necessary,
we may assume that
\[
  c_1=1,
  \qquad
  c_2=a>1.
\]
Set
\[
  \omega_1(t,r)=1+|r-t|,
  \qquad
  \omega_2(t,r)=1+|r-at|.
\]
Let \(X\) be the affine space of odd pairs \((U,V)\), with the prescribed
initial traces, such that
\[
  U,V\in C^2([0,\infty)\times\mathbb R)
\]
and
\[
\begin{aligned}
  \|(U,V)\|_X
  :={}&
  \left\|
    \omega_1^{\alpha_*}
    \bigl(
      |U_t|+|U_{tt}|+|U_{tr}|
    \bigr)
  \right\|_{L^\infty_{t,r}}
\\
  &+
  \left\|
    \omega_2^{\beta_*}
    \bigl(
      |V_t|+|V_{tt}|+|V_{tr}|
    \bigr)
  \right\|_{L^\infty_{t,r}}
  <\infty.
\end{aligned}
\]
Although the odd extensions are defined on the whole line, all weighted
\(L^\infty_{t,r}\)-norms in this proof are taken over
\((t,r)\in[0,\infty)\times[0,\infty)\). The odd extensions are used only
in the one-dimensional Duhamel formulas.
Define the nonlinear map
\[
  \mathcal S(U,V)
  =
  \bigl(
    U^0+L_1(r^{1-p}|V_t|^p),
    V^0+L_2(r^{1-q}|U_t|^q)
  \bigr).
\]
We first estimate the homogeneous and nonlinear parts of
\(\mathcal S\), and then prove the contraction property.

\medskip
\noindent\textbf{Step 1. Homogeneous estimates.}
\medskip

We claim that
\begin{equation}
\label{linear-estimate-global}
\begin{aligned}
  &\left\|
    \omega_1^{\alpha_*}
    \bigl(
      |U_t^0|+|U_{tt}^0|+|U_{tr}^0|
    \bigr)
  \right\|_{L^\infty_{t,r}}
+
  \left\|
    \omega_2^{\beta_*}
    \bigl(
      |V_t^0|+|V_{tt}^0|+|V_{tr}^0|
    \bigr)
  \right\|_{L^\infty_{t,r}}
  \leq C_0\vep
\end{aligned}
\end{equation}
for some constant \(C_0>0\).

Let \(c>0\), \(\mu>0\), and consider
\[
  Z_{tt}^0-c^2Z_{rr}^0=0
\]
with odd initial data
\[
  Z^0(0,r)=\phi_0(r),
  \qquad
  Z_t^0(0,r)=\phi_1(r),
\]
satisfying
\[
  |\phi_0'(r)|
  +|\phi_0''(r)|
  +|\phi_1(r)|
  +|\phi_1'(r)|
  \lesssim
  \vep(1+|r|)^{-\mu}.
\]
The d'Alembert formula gives
\[
  Z^0(t,r)
  =
  \frac12\bigl(\phi_0(r+ct)+\phi_0(r-ct)\bigr)
  +
  \frac1{2c}
  \int_{r-ct}^{r+ct}\phi_1(\lambda)\,d\lambda.
\]
Hence
\begin{align*}
  Z_t^0(t,r)
  &=
  \frac{c}{2}
  \bigl(
    \phi_0'(r+ct)-\phi_0'(r-ct)
  \bigr)
  +
  \frac12
  \bigl(
    \phi_1(r+ct)+\phi_1(r-ct)
  \bigr),
\\
  Z_{tt}^0(t,r)
  &=
  \frac{c^2}{2}
  \bigl(
    \phi_0''(r+ct)+\phi_0''(r-ct)
  \bigr)
  +
  \frac{c}{2}
  \bigl(
    \phi_1'(r+ct)-\phi_1'(r-ct)
  \bigr),
\\
  Z_{tr}^0(t,r)
  &=
  \frac{c}{2}
  \bigl(
    \phi_0''(r+ct)-\phi_0''(r-ct)
  \bigr)
  +
  \frac12
  \bigl(
    \phi_1'(r+ct)+\phi_1'(r-ct)
  \bigr).
\end{align*}
Since
\(
  |r+ct|\geq |r-ct|, r,t\geq0,
\)
we obtain
\[
  |Z_t^0(t,r)|
  +
  |Z_{tt}^0(t,r)|
  +
  |Z_{tr}^0(t,r)|
  \lesssim
  \vep(1+|r-ct|)^{-\mu}.
\]
By \eqref{ini}, applying the above estimate to
\[
  (U^0,c,\mu)=(U^0,1,\alpha_*),
  \qquad
  (V^0,c,\mu)=(V^0,a,\beta_*)
\]
yields \eqref{linear-estimate-global}.

\medskip
\noindent\textbf{Step 2. Nonlinear estimates.}
\medskip

Define
\[
  X_\vep
  =
  \left\{
    (U,V)\in X:
    \|(U,V)\|_X\leq 2C_0\vep
  \right\}.
\]
Let \((U,V)\in X_\vep\). We shall prove
\[
\begin{aligned}
  &\left\|
    \omega_1^{\alpha_*}
    \Bigl(
      \bigl|\partial_tL_1(r^{1-p}|V_t|^p)\bigr|
      +
      \bigl|\partial_{tt}L_1(r^{1-p}|V_t|^p)\bigr|
      +
      \bigl|\partial_{tr}L_1(r^{1-p}|V_t|^p)\bigr|
    \Bigr)
  \right\|_{L^\infty_{t,r}}
  \les
  \vep^p,
\\
  &\left\|
    \omega_2^{\beta_*}
    \Bigl(
      \bigl|\partial_tL_2(r^{1-q}|U_t|^q)\bigr|
      +
      \bigl|\partial_{tt}L_2(r^{1-q}|U_t|^q)\bigr|
      +
      \bigl|\partial_{tr}L_2(r^{1-q}|U_t|^q)\bigr|
    \Bigr)
  \right\|_{L^\infty_{t,r}}
  \les
  \vep^q.
\end{aligned}
\]
It suffices to prove the first estimate; the second follows by the same argument after interchanging the two components.

From \eqref{global-index} and \eqref{eq:scaling-indices}, we have
\[
  0<\alpha_*<p-1,
  \qquad
  0<\beta_*<q-1,
  \qquad
  0<p\beta_*<1,
  \qquad
  0<q\alpha_*<1,
\]
and
\[
  p\beta_*-\alpha_*=2-p>0,
  \qquad
  q\alpha_*-\beta_*=2-q>0.
\]
Set
\[
  F(s,r)=r^{1-p}|V_t(s,r)|^p\ .
\]
For \(r>0\),
\[
  F_r(s,r)
  =
  (1-p)r^{-p}|V_t(s,r)|^p
  +
  p r^{1-p}|V_t(s,r)|^{p-2}
  V_t(s,r)V_{tr}(s,r).
\]
Since \((U,V)\in X_\vep\),
\[
  |V_t(s,r)|+|V_{tr}(s,r)|
  \les
  \vep(1+|r-as|)^{-\beta_*}.
\]
Moreover, \(V_t(s,0)=0\), and hence, for \(0\leq r\leq1\),
\[
  |V_t(s,r)|
  \leq
  \int_0^r|V_{tr}(s,\rho)|\,d\rho
  \les
  \vep r(1+|r-as|)^{-\beta_*}.
\]
It follows that
\begin{equation}
\label{eq:F-Fr-bound}
  |F(s,r)|+|F_r(s,r)|
  \les
  \vep^p
  (1+r)^{1-p}
  (1+|r-as|)^{-p\beta_*},
  \qquad r\geq0.
\end{equation}
Since \(F\) is odd, \(F_r\) is even.
Differentiating the one-dimensional Duhamel formula \eqref{Duhamel operator} gives
\begin{align*}
  \partial_tL_1(F)(t,r)
  &=
  \frac12
  \int_0^t
  \Bigl[
    F(s,r+t-s)+F(s,r-t+s)
  \Bigr]\,ds,
\\
  \partial_{tr}L_1(F)(t,r)
  &=
  \frac12
  \int_0^t
  \Bigl[
    F_r(s,r+t-s)+F_r(s,r-t+s)
  \Bigr]\,ds,
\\
  \partial_{tt}L_1(F)(t,r)
  &=
  F(t,r)
  +
  \frac12
  \int_0^t
  \Bigl[
    F_r(s,r+t-s)-F_r(s,r-t+s)
  \Bigr]\,ds.
\end{align*}
By \eqref{eq:F-Fr-bound}, all the integral terms are bounded by the
same type of convolution. Applying Lemma~\ref{lem:oned-convolution}
with
\[
  \beta=p-1,
  \qquad
  \gamma=p\beta_*,
\]
and using
\[
  (p-1)+p\beta_*=1+\alpha_*,
\]
we obtain
\[
  |\partial_tL_1(F)(t,r)|
  +
  |\partial_{tr}L_1(F)(t,r)|
  \les
  \vep^p(1+|r-t|)^{-\alpha_*},
\]
and the same bound holds for the integral term in
\(\partial_{tt}L_1(F)\). Moreover, by \eqref{eq:F-Fr-bound},
\[
  |F(t,r)|
  \les
  \vep^p
  (1+r)^{-(p-1)}
  (1+|r-at|)^{-p\beta_*}
  \les
  \vep^p(1+|r-t|)^{-\alpha_*},
\]
where we used
\[
  \alpha_*<\min\{p-1,p\beta_*\},
  \qquad
  1+|r-t|
  \les
  \max\{1+r,1+|r-at|\}.
\]
Hence
\[
  |\partial_tL_1(F)|
  +
  |\partial_{tt}L_1(F)|
  +
  |\partial_{tr}L_1(F)|
  \les
  \vep^p(1+|r-t|)^{-\alpha_*}.
\]
The estimate for the second component follows in the same way, with
\[
  p,\ \beta_*,\ V
  \qquad\text{replaced by}\qquad
  q,\ \alpha_*,\ U,
\]
and with the change of variables \(\tau=as\) in
Lemma~\ref{lem:oned-convolution}. Here we use
\[
  (q-1)+q\alpha_*=1+\beta_*,
  \qquad
  \beta_*<q-1,
  \qquad
  \beta_*<q\alpha_*.
\]
Consequently,
\begin{equation}
\label{eq:L2-uniform}
\begin{aligned}
  &\left\|
    \omega_1^{\alpha_*}
    \Bigl(
      \bigl|\partial_tL_1(r^{1-p}|V_t|^p)\bigr|
      +
      \bigl|\partial_{tt}L_1(r^{1-p}|V_t|^p)\bigr|
      +
      \bigl|\partial_{tr}L_1(r^{1-p}|V_t|^p)\bigr|
    \Bigr)
  \right\|_{L^\infty_{t,r}}
\\
  \qquad&+
  \left\|
    \omega_2^{\beta_*}
    \Bigl(
      \bigl|\partial_tL_2(r^{1-q}|U_t|^q)\bigr|
      +
      \bigl|\partial_{tt}L_2(r^{1-q}|U_t|^q)\bigr|
      +
      \bigl|\partial_{tr}L_2(r^{1-q}|U_t|^q)\bigr|
    \Bigr)
  \right\|_{L^\infty_{t,r}}
  \les
  \vep^p+\vep^q.
\end{aligned}
\end{equation}
Combining \eqref{linear-estimate-global} and
\eqref{eq:L2-uniform}, we obtain
\beeq\label{eq:L2-uniform-2}
  \|\mathcal S(U,V)\|_X
  \leq
  C_0\vep+C(\vep^p+\vep^q).
\eneq

\medskip
\noindent\textbf{Step 3. Persistence of regularity.}
\medskip

It remains to verify that the map \(\mathcal S\) preserves the
\(C^2\)-regularity required in the definition of \(X\). Let
\((U,V)\in X_\vep\). Since \(V_t(t,\cdot)\) is odd and
\(V_t\in C^1\), we have
\[
  \frac{V_t(t,r)}{r}
  =
  \int_0^1 V_{tr}(t,\theta r)\,d\theta
  \longrightarrow
  V_{tr}(t,0)
  \qquad\text{as }r\to0.
\]
Hence
\[
  \frac{V_t}{r}\in C([0,\infty)\times\mathbb R).
\]
Write
\[
  F(t,r)
  =
  r\left|\frac{V_t(t,r)}{r}\right|^p\ .
\]
For \(r\neq0\),
\[
\begin{aligned}
  F_r
  =
  (1-p)
  \left|\frac{V_t}{r}\right|^p
  +
  p
  \left|\frac{V_t}{r}\right|^{p-2}
  \frac{V_t}{r}V_{tr}\ ,\quad
  F_t
  =
  p
  \left|\frac{V_t}{r}\right|^{p-2}
  \frac{V_t}{r}V_{tt}.
\end{aligned}
\]
Since \(V_{tt}(t,\cdot)\) is odd,
\(
  V_{tt}(t,0)=0.
\)
Therefore,
\[
  F_r(t,r)\longrightarrow |V_{tr}(t,0)|^p,
  \qquad
  F_t(t,r)\longrightarrow0
  \qquad\text{as }r\to0.
\]
Defining
\[
  F_r(t,0)=|V_{tr}(t,0)|^p,
  \qquad
  F_t(t,0)=0,
\]
we obtain
\[
  F_r,F_t\in C([0,\infty)\times\mathbb R),
\]
and hence
\[
  F\in C^1([0,\infty)\times\mathbb R).
\]
The same argument gives
\[
  r^{1-q}|U_t|^q
  \in C^1([0,\infty)\times\mathbb R).
\]
Therefore,  
\[
  L_1(r^{1-p}|V_t|^p),
  \quad
  L_2(r^{1-q}|U_t|^q)
  \in C^2([0,\infty)\times\mathbb R).
\]
Consequently,
\[
  \mathcal S(U,V)\in X.
\]

Combining this regularity property with \eqref{eq:L2-uniform-2}, we obtain,
for \(\vep>0\) sufficiently small,
\[
  \|\mathcal S(U,V)\|_X
  \leq
  C_0\vep+C(\vep^p+\vep^q)
  \leq
  2C_0\vep.
\]
Thus
\[
  \mathcal S(X_\vep)\subset X_\vep.
\]

\medskip
\noindent\textbf{Step 4. Contraction.}
\medskip

We equip \(X_\vep\) with the weaker distance
\[
\begin{aligned}
  d\bigl((U_1,V_1),(U_2,V_2)\bigr)
  :={}&
  \left\|
    \omega_1^{\alpha_*}
    \bigl((U_1)_t-(U_2)_t\bigr)
  \right\|_{L^\infty_{t,r}}
+
  \left\|
    \omega_2^{\beta_*}
    \bigl((V_1)_t-(V_2)_t\bigr)
  \right\|_{L^\infty_{t,r}}.
\end{aligned}
\]
Let \((U_i,V_i)\in X_\vep\), \(i=1,2\). Using
\[
  \bigl||x|^p-|y|^p\bigr|
  \les
  (|x|+|y|)^{p-1}|x-y|,
\]
together with the bounds defining \(X_\vep\) and
\((V_i)_t(s,0)=0\), we obtain
\[
\begin{aligned}
 &r^{1-p}
 \bigl|
   |(V_1)_t|^p-|(V_2)_t|^p
 \bigr|
\les
 \vep^{p-1}
 (1+r)^{1-p}
 \omega_2^{-p\beta_*}
 d\bigl((U_1,V_1),(U_2,V_2)\bigr).
\end{aligned}
\]
Hence, by the same convolution estimate as in Step~2,
\[
\begin{aligned}
 &\left\|
   \omega_1^{\alpha_*}
   \partial_tL_1
   \left(
     r^{1-p}
     \bigl(
       |(V_1)_t|^p-|(V_2)_t|^p
     \bigr)
   \right)
 \right\|_{L^\infty_{t,r}}
\les
 \vep^{p-1}
 d\bigl((U_1,V_1),(U_2,V_2)\bigr).
\end{aligned}
\]
The corresponding estimate for the second component gives
\[
\begin{aligned}
 &\left\|
   \omega_2^{\beta_*}
   \partial_tL_2
   \left(
     r^{1-q}
     \bigl(
       |(U_1)_t|^q-|(U_2)_t|^q
     \bigr)
   \right)
 \right\|_{L^\infty_{t,r}}
\les
 \vep^{q-1}
 d\bigl((U_1,V_1),(U_2,V_2)\bigr).
\end{aligned}
\]
Therefore,
\[
  d\bigl(
    \mathcal S(U_1,V_1),
    \mathcal S(U_2,V_2)
  \bigr)
  \les
  \bigl(
    \vep^{p-1}+\vep^{q-1}
  \bigr)
  d\bigl(
    (U_1,V_1),(U_2,V_2)
  \bigr).
\]
Fix \(\vep>0\) sufficiently small that the last contraction factor is
strictly less than one. 

Choose \((U_0,V_0)\in X_\vep\) and define
\[
  (U_{n+1},V_{n+1})
  =
  \mathcal S(U_n,V_n),
  \qquad n\geq0.
\]
Then, the contraction estimate gives
\[
  (U_n)_t\longrightarrow U_t,
  \qquad
  (V_n)_t\longrightarrow V_t\ , \qquad n\to\infty\ ,
\]
in the corresponding weighted \(L^\infty\)-norms. Since all iterates
have the same prescribed initial traces, it follows that
\[
  U_n\longrightarrow U,
  \qquad
  V_n\longrightarrow V\ , \qquad n\to \infty\ ,
\]
locally uniformly on \([0,\infty)\times\mathbb R\). The difference
estimates above imply, locally uniformly,
\[
  r^{1-p}|(V_n)_t|^p
  \longrightarrow
  r^{1-p}|V_t|^p,
  \qquad
  r^{1-q}|(U_n)_t|^q
  \longrightarrow
  r^{1-q}|U_t|^q.
\]
Hence, passing to the limit in the Duhamel formulas, we obtain
\[
  (U,V)=\mathcal S(U,V).
\]

It remains to recover the \(C^2\)-regularity. Let
\((U^{\mathrm{cl}},V^{\mathrm{cl}})\) be the local classical solution
with the same initial data, defined on its maximal interval
\([0,T_{\mathrm{cl}})\). By the local version of the preceding
difference estimate, the integral system is unique in the class of
solutions with bounded weighted time derivatives and linear vanishing
at \(r=0\). Therefore,
\[
  (U^{\mathrm{cl}},V^{\mathrm{cl}})
  =
  (U,V)
  \qquad\text{on }[0,T_{\mathrm{cl}}).
\]
The uniform estimates defining \(X_\vep\) give uniform local
Lipschitz bounds for \((U_n)_t\) and \((V_n)_t\). Passing to the
locally uniform limit shows that \(U_t\) and \(V_t\) are locally
Lipschitz. Since the local classical solution agrees with \((U,V)\),
these bounds hold classically on \([0,T_{\mathrm{cl}})\). Hence the
standard continuation criterion for the differentiated integral
system excludes a finite \(T_{\mathrm{cl}}\). Consequently,
\[
  T_{\mathrm{cl}}=\infty,
  \qquad
  U,V\in C^2([0,\infty)\times\mathbb R).
\]
Since
\[
  u(t,r)
  =
  \frac{U(t,r)}{r}
  =
  \int_0^1 U_r(t,\theta r)\,d\theta,
  \qquad
  v(t,r)
  =
  \frac{V(t,r)}{r}
  =
  \int_0^1 V_r(t,\theta r)\,d\theta,
\]
we have
\[
  u_t(t,r)
  =
  \int_0^1 U_{tr}(t,\theta r)\,d\theta,
  \qquad
  v_t(t,r)
  =
  \int_0^1 V_{tr}(t,\theta r)\,d\theta,
\]
and
\[
  u_r(t,r)
  =
  \int_0^1 \theta U_{rr}(t,\theta r)\,d\theta,
  \qquad
  v_r(t,r)
  =
  \int_0^1 \theta V_{rr}(t,\theta r)\,d\theta.
\]
Thus \(u\) and \(v\) extend across \(r=0\) as \(C^1\)-functions, with
\[
  u_r(t,0)=v_r(t,0)=0.
\]
Therefore, \((u,v)\) is a global radial \(C^1\)-solution of
\eqref{eq:main-system}, completing the proof.
\end{proof}

\section{Finite-Time Blow-Up}
\label{sec:finite-time-blow-up}

We now turn to the finite-time blow-up part of Theorem~\ref{thm:main}.
The proof begins with the construction of suitable functionals for the two
components. The required lower bounds are obtained from the contribution
of the initial data along the corresponding outgoing characteristic cones
and are then combined with the differential inequalities satisfied by these
functionals. In the subcritical case, a standard Kato-type argument suffices,
while the critical case requires an additional logarithmic improvement.
We therefore first recall the subcritical Kato-type lemma and then establish
a critical version adapted to our argument.

\subsection{Subcritical and critical Kato-type lemmas}
\label{subsec:blowup-kato}

\begin{lemma}[Subcritical power-type Kato lemma]
\label{lem:subcritical-power-Kato}
Let \(\nu>1\), \(b\geq1\), and \(\mu\in\mathbb R\), with
\[
 D:=(\nu-1)b-\mu+2>0.
\]
Suppose that \(H\in C^2([0,T))\) is nonnegative and satisfies
\[
\begin{cases}
H(0) \geq 0,
  &H'(0) \geq 0,
  \\
  H''(t) \geq C_0(1+t)^{-\mu}H(t)^\nu,
  &0\leq t<T,
  \\
  H(t) \geq C_1 A t^b,
  &T_0\leq t<T
\end{cases},
\]
where \(C_0,C_1>0\), \(T_0\geq1\), and \(A>0\).
There exist \(A_0>0\) and \(C>0\), depending only on
\(C_0,C_1,T_0,\mu,\nu,b\), such that
\begin{equation*}
 T\leq C A^{-\frac{\nu-1}{D}},\qquad 0<A\leq A_0.
\end{equation*}
In particular, no such function can exist for all \(t\geq0\).
\end{lemma}
The subcritical case \(D>0\) is classical; see Lemma~4.1 in
\cite{L-W} for the above formulation. We also refer to
Sideris~\cite{S-84} for the blow-up argument and
Takamura~\cite{Taka15} for the corresponding upper bound for the lifespan.

At the critical threshold
\[
  D=(\nu-1)b-\mu+2=0,
\]
the above power-type estimate is no longer applicable. For \(b=1\),
a quantitative upper bound for the lifespan was obtained in
\cite{Liu2026Kato}. A critical Kato-type estimate for \(b\geq1\) was
also established in \cite{LPT26}, under an additional largeness
condition on the constant in the lower bound. In the present argument,
we instead use a logarithmically improved lower bound of the form
\[
  H(t)\geq C_1At^b\log(1+t),
\]
which leads to the following critical version.

\begin{lemma}[Critical logarithmic Kato lemma]
\label{lem:critical-log-Kato}
Let \(\nu>1\), \(b>0\), and
\[
 \mu=2+b(\nu-1).
\]
Suppose that \(H\in C^2([0,T))\) satisfies, for \(T_0\leq t<T\),
\[
\begin{cases}
  H(t) > 0,
  \quad H'(t) \geq 0,
  \\
  H''(t) \geq C_0(1+t)^{-\mu}H(t)^\nu,
  \\
  H(t) \geq C_1 A t^b\log(1+t)
\end{cases},
\]
where \(C_0,C_1>0\) and \(T_0\geq1\).
There exist \(A_0>0\) and \(C>0\), depending only on
\(C_0,C_1,T_0,b,\nu\), such that
\begin{equation}
\label{eq:critical-log-Kato-lifespan}
 T\leq\exp(C/A),\qquad 0<A\leq A_0.
\end{equation}
\end{lemma}

\begin{proof}
Let \(R\geq T_0\) be chosen below, and suppose that \(T>2R\).
Define
\[
 Y(s)=R^{-b}H(Rs),\qquad 1\leq s\leq2.
\]
Then \(Y'\geq0\), and
\[
 Y''(s)\geq
 C_0R^{2-b+b\nu}(1+Rs)^{-\mu}Y(s)^\nu.
\]
Since \(1+Rs\leq2Rs\), \(s\leq2\), and
\(2-b-\mu+b\nu=0\), we obtain
\begin{equation}
\label{eq:Y-autonomous-lower}
 Y''(s)\geq C_2Y(s)^\nu,\qquad
 C_2=C_0 2^{-2\mu}>0.
\end{equation}
Moreover,
\begin{equation}
\label{eq:Y-initial-lower}
 Y(1)=R^{-b}H(R)\geq C_1 A\log(1+R).
\end{equation}
Multiplying \eqref{eq:Y-autonomous-lower} by \(Y'\) and integrating,
we find
\[
 Y'(s)\geq
 \left(\frac{2C_2}{\nu+1}\right)^{1/2}
 \left(Y(s)^{\nu+1}-Y(1)^{\nu+1}\right)^{1/2},
 \qquad s>1.
\]
Consequently,
\begin{align}
 s-1
 &\leq
 \left(\frac{\nu+1}{2C_2}\right)^{1/2}
 \int_{Y(1)}^{Y(s)}
 \frac{d\xi}{\left(\xi^{\nu+1}-Y(1)^{\nu+1}\right)^{1/2}}
 \notag\\
 &\leq C_3Y(1)^{-(\nu-1)/2},
 \label{eq:Y-time-upper}
\end{align}
where
\[
 C_3=
 \left(\frac{\nu+1}{2C_2}\right)^{1/2}
 \int_1^\infty\frac{dz}{(z^{\nu+1}-1)^{1/2}}<\infty.
\]
The integrand has order \((z-1)^{-1/2}\) near \(1\) and
\(z^{-(\nu+1)/2}\) at infinity, so this improper integral converges.

Now choose
\[
 R=\exp\left(\frac{C_3^{2/(\nu-1)}}{C_1 A}\right).
\]
For \(A_0\) sufficiently small, \(R\geq T_0\). By
\eqref{eq:Y-initial-lower},
\[
 Y(1)\geq C_1 A\log(1+R)>C_3^{2/(\nu-1)}.
\]
Taking \(s=2\) in \eqref{eq:Y-time-upper} gives
\[
 1\leq C_3Y(1)^{-(\nu-1)/2}<1,
\]
a contradiction. Thus \(T\leq2R\), and increasing \(C\) gives
\eqref{eq:critical-log-Kato-lifespan}.
\end{proof}

In the following, let \((u,v)\) be the corresponding radial
\(C^1\)-solution of \eqref{eq:main-system}. We assume that
\[
  f_1=f_2=0,
  \qquad
  \supp g_i\subset B_R,
  \qquad
  g_i\geq0,\quad g_i\not\equiv0,
  \qquad i=1,2.
\]
By positivity of the radial representation and the nonnegativity of
the source terms, we have
\[
 u(t,r)\geq0,
 \qquad
 v(t,r)\geq0,
 \qquad t\geq0,\quad r\geq0.
\]

After relabeling the components and rescaling space by the larger
propagation speed if necessary, we normalize
\[
  c_1=c\in(0,1),
  \qquad
  c_2=1.
\]
By finite speed of propagation, the solution satisfies
\[
  \supp(u,v,u_t,v_t)(t,\cdot)
  \subset B_{R+t}.
\]

\subsection{Spatial integrals and subcritical blow-up}
\label{subsec:subcritical-blowup}

We define the functionals
\begin{equation*}
 F(t)=\int_{\mathbb R^3}u_t(t,x)\,dx,\qquad
 G(t)=\int_{\mathbb R^3}v_t(t,x)\,dx.
\end{equation*}
Integrating both equations in \eqref{eq:main-system} over $\R^3$ yields 
\begin{equation}
\label{eq:FG-exact}
 F'(t)=\int_{\mathbb R^3}|v_t(t,x)|^p\,dx,\qquad
 G'(t)=\int_{\mathbb R^3}|u_t(t,x)|^q\,dx.
\end{equation}
For the above radial solutions, these identities hold in the absolutely
continuous sense. The nonnegative, nontrivial initial velocities imply
\[
 F(0)=\vep\int_{\mathbb R^3}g_1\,dx>0,\qquad
 G(0)=\vep\int_{\mathbb R^3}g_2\,dx>0.
\]
Thus \(F,G\) are positive and nondecreasing.

By H\"older's inequality and \eqref{eq:FG-exact}, we have 
\[
 G(t)^p
 \leq F'(t)|B_{R+t}|^{p-1}
 \lesssim (1+t)^{3(p-1)}F'(t),
\]
and the analogous estimate with the components interchanged. Hence
\begin{equation}
\label{eq:coupled-ode}
\left\{
\begin{aligned}
 F'(t)&\gtrsim(1+t)^{-3(p-1)}G(t)^p,\\
 G'(t)&\gtrsim(1+t)^{-3(q-1)}F(t)^q.
\end{aligned}
\right.
\end{equation}

We next establish preliminary power-type lower bounds. Since \(g_2\geq0\)
and \(g_2\not\equiv0\), there are
\(0<\rho_0<\rho_1<R\) and \(d_0>0\) such that
\begin{equation*}
 rg_2(r)\geq d_0,\qquad \rho_0\leq r\leq\rho_1.
\end{equation*}
Therefore, for \(r-c_2t\in [\rho_0, \rho_1]\), 
\begin{equation}
\label{eq:Vt-shell}
 V_t(t,r)\geq V_t^0(t,r)=\frac{\vep}{2}
 \left[(r+c_2t)g_2(r+c_2t)+(r-c_2t)g_2(r-c_2t)\right]
  \geq\frac{\ep d_0}{2}.
\end{equation}
Since \(v_t=V_t/r\), it follows that
\[
 F'(t)=4\pi\int_0^\infty r^{2-p}|V_t(t,r)|^p\,dr
 \gtrsim\vep^p
 \int_{c_2t+\rho_0}^{c_2t+\rho_1}r^{2-p}\,dr
 \gtrsim\vep^p(1+t)^{2-p}.
\]
Integrating over \([t/2,t]\), and arguing similarly for \(U_t\),
gives, for \(t\geq2\),
\begin{equation*}
 F(t)\gtrsim\vep^p t^{3-p},\qquad
 G(t)\gtrsim\vep^q t^{3-q}.
\end{equation*}
To transform the system \eqref{eq:coupled-ode} into second-order differential inequalities, we set
\[
 H(t)=\int_0^t F(s)\,ds,\qquad
 K(t)=\int_0^t G(s)\,ds.
\]
Integrating the second inequality in \eqref{eq:coupled-ode} and
applying H\"older's inequality in time, we obtain
\begin{align*}
 G(t)
 &\gtrsim(1+t)^{-3(q-1)}\int_0^t F(s)^q\,ds
 \gtrsim(1+t)^{-4(q-1)}H(t)^q.
\end{align*}
Substitution into the
first inequality in \eqref{eq:coupled-ode} gives
\begin{equation}
\label{eq:H-scalar-ode}
 H''(t)\gtrsim(1+t)^{-(4pq-p-3)}H(t)^{pq}.
\end{equation}
By symmetry,
\begin{equation}
\label{eq:K-scalar-ode}
 K''(t)\gtrsim(1+t)^{-(4pq-q-3)}K(t)^{pq}.
\end{equation}
The initial lower bounds imply that, for \(t\geq t_0>1\),
\begin{equation*}
 H(t)\gtrsim\vep^p t^{4-p},\qquad
 K(t)\gtrsim\vep^q t^{4-q}.
\end{equation*}

Suppose first that \(p\leq q\). For \eqref{eq:H-scalar-ode}, set
\[
 \nu=pq,\qquad \mu=4pq-p-3,\qquad b=4-p,\qquad A=\vep^p.
\]
Then the parameter in Lemma~\ref{lem:subcritical-power-Kato} is
\[
 D=(pq-1)(4-p)-(4pq-p-3)+2
  =p\left(2+\frac1p-pq\right)=p\Gamma(p,q)>0.
\]
Applying that lemma gives
\[
 T_\vep\lesssim
 (\vep^p)^{-\frac{pq-1}{p\Gamma(p,q)}}
 =\vep^{-\frac{pq-1}{\Gamma(p,q)}}.
\]
If \(q\leq p\), we instead use \eqref{eq:K-scalar-ode} with
\[
 \nu=pq,\qquad \mu=4pq-q-3,\qquad b=4-q,\qquad A=\vep^q.
\]
In this case,
\[
 D=(pq-1)(4-q)-(4pq-q-3)+2
  =q\left(2+\frac1q-pq\right)=q\Gamma(p,q)>0.
\]
Lemma~\ref{lem:subcritical-power-Kato} therefore yields
\[
 T_\vep\lesssim
 (\vep^q)^{-\frac{pq-1}{q\Gamma(p,q)}}
 =\vep^{-\frac{pq-1}{\Gamma(p,q)}}.
\]
This proves the first inequality in \eqref{eq:compact-blowup-lifespan}.

\subsection{Critical blow-up}
\label{subsec:critical-blowup}

In this subsection, we prove finite-time blow-up in the critical case.
The key step is to improve the lower bounds for \(F\) and \(G\) on the
critical curve by an additional logarithmic factor. Once these improved
lower bounds are established, the upper bound for the lifespan in the critical case follows
from Lemma~\ref{lem:critical-log-Kato}.

We first state the logarithmic improvement.

\begin{lemma}[Logarithmic improvement on the critical curve]
\label{lem:critical-log-improvement}
There exist constants \(C>0\) and \(t_0\geq1\), independent of \(\vep\),
such that
\begin{equation}
\label{eq:critical-log-lower-F}
  F(t)
  \geq
  C\vep^{p^2q}t^{3-p}\log(1+t),
  \qquad
  p\leq q,
  \qquad
  pq=2+\frac1p,
\end{equation}
and
\begin{equation}
\label{eq:critical-log-lower-G}
  G(t)
  \geq
  C\vep^{pq^2}t^{3-q}\log(1+t),
  \qquad
  q<p,
  \qquad
  pq=2+\frac1q,
\end{equation}
for \(t_0\leq t<T_\vep\).
\end{lemma}
With the help of Lemma~\ref{lem:critical-log-improvement}, we can now
derive the upper bound for the lifespan in the critical case.
If \(p\leq q\) and \(pq=2+1/p\), then
\eqref{eq:critical-log-lower-F} yields
\[
  H(t)
  \gtrsim
  \vep^{p^2q}t^{4-p}\log(1+t).
\]
Applying Lemma~\ref{lem:critical-log-Kato} to
\eqref{eq:H-scalar-ode} with
\[
  \nu=pq,\qquad
  b=4-p,\qquad
  \mu=4pq-p-3,\qquad
  A=\vep^{p^2q},
\]
and noting that
\[
  \mu-\bigl(2+b(\nu-1)\bigr)
  =
  p\left(pq-2-\frac1p\right)=0,
\]
we obtain
\begin{equation*}
  T_\vep
  \leq
  \exp\left(C\vep^{-p^2q}\right)
  =
  \exp\left(C\vep^{-(2p+1)}\right).
\end{equation*}
If \(q<p\) and \(pq=2+1/q\), then
\eqref{eq:critical-log-lower-G} gives
\[
  K(t)
  \gtrsim
  \vep^{pq^2}t^{4-q}\log(1+t).
\]
Applying Lemma~\ref{lem:critical-log-Kato} to
\eqref{eq:K-scalar-ode} with
\[
  \nu=pq,\qquad
  b=4-q,\qquad
  \mu=4pq-q-3,\qquad
  A=\vep^{pq^2},
\]
and using
\[
  \mu-\bigl(2+b(\nu-1)\bigr)
  =
  q\left(pq-2-\frac1q\right)=0,
\]
we obtain
\begin{equation*}
  T_\vep
  \leq
  \exp\left(C\vep^{-pq^2}\right).
\end{equation*}
Thus \eqref{eq:compact-blowup-lifespan} follows.

It remains to prove Lemma~\ref{lem:critical-log-improvement}.
On the critical curve, the power lower bounds for \(F\) and \(G\)
make the parameter in the Kato-type argument vanish. We therefore
seek an additional logarithmic factor. When \(p\leq q\), we use a
superposition of exponential eigenfunctions to obtain a logarithmic
lower bound for a weighted pairing with \(v_t\). When \(q<p\), we use
a test function transported along the slower characteristic cone.
The faster-wave Duhamel region intersects a slow outgoing shell over
a time interval of length comparable to \(t\); the resulting source
lower bound, integrated over spatial scales, yields a logarithmic
lower bound for a weighted pairing with \(u_t\). In each case,
H\"older's inequality converts the pairing estimate into the
logarithmically improved lower bound required by
Lemma~\ref{lem:critical-log-Kato}.

\subsubsection{\textbf{Logarithmic improvement of the lower bound for \(F\)}}

The logarithmic improvement of the lower bound for \(F\) is established under the assumptions
\beeq\label{critical-p}
  p\leq q,
  \qquad
  pq=2+\frac1p.
\eneq
In particular, \(1<p<2\).  For $\la\in (0, 1)$, let
\[
 \phi_\lambda(x)=\frac{\sinh(\lambda|x|)}{\lambda|x|},
 \qquad \phi_\lambda(0)=1,\qquad
 \Phi_\lambda(t,x)=e^{-\lambda t}\phi_\lambda(x) .
\]
Then \(\Delta\phi_\lambda=\lambda^2\phi_\lambda\).
Define
\begin{equation}
\label{eq:Phi-critical-p}
 \Phi(t,x)=\int_0^1\lambda^{1/p'}\Phi_\lambda(t,x)\,d\lambda.
\end{equation}

We first establish a logarithmic lower bound for the pairing of \(v_t\)
with \(\Phi\).

\begin{lemma}
\label{lem:critical-pairing-p}
Assume \eqref{critical-p}. Then there exist
\(t_1\geq1\) and \(C>0\), independent of \(\vep\), such that
\begin{equation*}
 \int_{\mathbb R^3}v_t(t,x)\Phi(t,x)\,dx
 \geq C\vep^{pq}\log(1+t),
 \qquad t_1\leq t<T_\vep.
\end{equation*}
\end{lemma}

\begin{proof}
For \(0<\lambda\leq1\), we introduce the following auxiliary functions:
\begin{align*}
 H_{1,\lambda}(t)
 &=\int_{\mathbb R^3}u_t(t,x)\Phi_\lambda(t,x)\,dx,
 \\
 H_{2,\lambda}(t)
 &=\int_{\mathbb R^3}v_t(t,x)\Phi_\lambda(t,x)\,dx,
\\
 G_\lambda(t)
 &=\int_{\mathbb R^3}(v_t+\lambda v)(t,x)\Phi_\lambda(t,x)\,dx.
\end{align*}
Direct computation gives
\[
 (\partial_t+2\lambda)
 (2H_{2,\lambda}-G_\lambda)
 =\int_{\mathbb R^3}|u_t|^q\Phi_\lambda\,dx\geq0.
\]
Since \(f_i=0\), we have
\[
 (2H_{2,\lambda}-G_\lambda)(0)
 = \vep\int_{\mathbb R^3}g_2\Phi_\lambda\,dx>0.
\]
It follows that
\[
 H_{2,\lambda}(t)\geq\frac12G_\lambda(t)\geq0,
 \qquad 0<\lambda\leq1.
\]
Here we have used
\[
 G_\lambda'(t)
 =\int_{\mathbb R^3}|u_t|^q\Phi_\lambda\,dx\geq0.
\]
The first equation in \eqref{eq:main-system}, with \(c_1=c\), gives
\[
 H_{1,\lambda}'+\lambda H_{1,\lambda}
 =c^2\lambda^2\int_{\mathbb R^3}u\Phi_\lambda\,dx
  +\int_{\mathbb R^3}|v_t|^p\Phi_\lambda\,dx
 \geq\int_{\mathbb R^3}|v_t|^p\Phi_\lambda\,dx.
\]
In particular, \(H_{1,\lambda}\geq0\).
By finite speed of propagation and the explicit form of \(\phi_\lambda\),
\[
 \int_{|x|\leq t+R}\Phi_\lambda(t,x)\,dx
 \lesssim\lambda^{-2}(1+t),\qquad 0<\lambda\leq1.
\]
H\"older's inequality therefore yields
\[
 G_\lambda'(t)
 \gtrsim\lambda^{2(q-1)}(1+t)^{-(q-1)}
 H_{1,\lambda}(t)^q.
\]
Collecting the preceding estimates, we obtain
\begin{equation}
\label{eq:core-estimates-critical-p}
\left\{
\begin{aligned}
  &H_{1,\lambda}'+\lambda H_{1,\lambda}
  \gtrsim
  \lambda^{2(p-1)}(1+t)^{-(p-1)}
  H_{2,\lambda}^p,
\\
 & H_{2,\lambda}(t)
  \geq
  \frac12 G_\lambda(t),
  \quad
  G_\lambda'(t)
  \gtrsim
  \lambda^{2(q-1)}(1+t)^{-(q-1)}
  H_{1,\lambda}^q.
\end{aligned}
\right.
\end{equation}

We next use the outgoing-shell lower bound to initiate the above
coupled estimates. On the outgoing shell
\(t+\rho_0\leq r\leq t+\rho_1\),
\eqref{eq:Vt-shell} gives \(v_t(t,r)\gtrsim\vep t^{-1}\) for
\(t\geq1\). If \(\lambda^{-1}\leq t\leq2\lambda^{-1}\), then
\(\Phi_\lambda(t,r)\gtrsim1\) on this shell. Hence
\[
 H_{1,\lambda}'+\lambda H_{1,\lambda}
 \gtrsim\vep^p t^{2-p}
 \gtrsim\vep^p\lambda^{p-2}.
\]
Multiplying by \(e^{\lambda t}\) and integrating from
\(\lambda^{-1}\) to \(2\lambda^{-1}\), we obtain
\[
 H_{1,\lambda}(2\lambda^{-1})
 \gtrsim\vep^p\lambda^{p-3}.
\]
The first inequality in \eqref{eq:core-estimates-critical-p} then gives
\begin{equation*}
 H_{1,\lambda}(t)\gtrsim\vep^p\lambda^{p-3},
 \qquad 2\lambda^{-1}\leq t\leq3\lambda^{-1}.
\end{equation*}
On this interval,
\[
 G_\lambda'(t)
 \gtrsim\lambda^{3(q-1)}
          (\vep^p\lambda^{p-3})^q
 =\vep^{pq}\lambda^{pq-3}.
\]
Integrating over \([2\lambda^{-1},3\lambda^{-1}]\), and using
monotonicity thereafter, yields
\begin{equation*}
 H_{2,\lambda}(t)
 \gtrsim\vep^{pq}\lambda^{pq-4},
 \qquad t\geq3\lambda^{-1}.
\end{equation*}
Since
\[
 pq-4=-1-\frac1{p'},
\]
Fubini's theorem and \(H_{2,\lambda}\geq0\) for every
\(0<\lambda\leq1\) give, for \(t>t_1>3\),
\begin{align*}
 \int_{\mathbb R^3}v_t\Phi\,dx
 &=\int_0^1\lambda^{1/p'}H_{2,\lambda}(t)\,d\lambda\\
 &\geq\int_{3/t}^1\lambda^{1/p'}H_{2,\lambda}(t)\,d\lambda\\
 &\gtrsim\vep^{pq}\int_{3/t}^1\frac{d\lambda}{\lambda}
 \gtrsim\vep^{pq}\log(1+t).
\end{align*}
\end{proof}

To convert the preceding pairing estimate into a lower bound for
\(F\), we also need the following \(L^{p'}\) estimate for \(\Phi\).

\begin{lemma}[\(L^{p'}\) estimate for $\Phi$]
\label{lem:Lp-Phi}
Let \(\Phi\) be defined by
\eqref{eq:Phi-critical-p}. For every fixed \(R>0\),
\begin{equation*}
 \int_{|x|\leq t+R}\Phi(t,x)^{p'}\,dx
 \leq C(1+t)^{2-p'}\log(2+t),\qquad t\geq1,
\end{equation*}
where \(C\) depends only on \(p,R\).
\end{lemma}

\begin{proof}
Since
\(
 \phi_\lambda(r)\lesssim e^{\lambda r}/(1+\lambda r),
\)
we have
\begin{equation*}
 \Phi(t,r)\lesssim
 \int_0^1\frac{\lambda^{1/p'} e^{-\lambda(t-r)}}{1+\lambda r}
 \,d\lambda.
\end{equation*}
For \(0\leq r\leq t/2\),
\[
 \Phi(t,r)\lesssim
 \int_0^\infty\lambda^{1/p'} e^{-\lambda t/2}\,d\lambda
 \lesssim(1+t)^{-1-1/p'}.
\]
For \(t/2\leq r\leq t+R\), we have \(r\simeq1+t\),
\[
 \Phi(t,r)\lesssim
 r^{-1}
 \int_0^1\lambda^{1/p'-1}e^{-\lambda (t+R-r)}\,d\lambda
 \lesssim(1+t)^{-1}(1+t+R-r)^{-1/p'}.
\]
Hence,
\begin{align*}
&\int_{|x|\leq t+R}\Phi(t,x)^{p'}\,dx\\
\les &\int_0^{t/2}r^2(1+t)^{-p'-1}\,dr
+(1+t)^{2-p'}
 \int_{t/2}^{t+R}(1+t+R-r)^{-1}\,dr\\
\lesssim & \ (1+t)^{2-p'}\log(2+t).
\end{align*}
\end{proof}

We now combine the preceding two lemmas. By H\"older's inequality,
\[
 \left(\int_{\mathbb R^3}v_t\Phi\,dx\right)^p
 \leq F'(t)
 \left(\int_{|x|\leq t+R}\Phi^{p'}\,dx\right)^{p-1}.
\]
Lemmas~\ref{lem:critical-pairing-p} and~\ref{lem:Lp-Phi} give
\begin{align*}
 F'(t)
 \gtrsim
 \frac{\vep^{p^2q}(\log(1+t))^p}
      {((1+t)^{2-p'}\log(2+t))^{p-1}}
 \gtrsim\vep^{p^2q}(1+t)^{2-p}\log(1+t)\ .
\end{align*}
Integrating over \([t/2,t]\), we obtain
\[
  F(t)
  \gtrsim
  \vep^{p^2q}t^{3-p}\log(1+t)\ ,
\]
for all sufficiently large \(t\). This proves
\eqref{eq:critical-log-lower-F} in Lemma \ref{lem:critical-log-improvement}.

\subsubsection{\textbf{Logarithmic improvement of the lower bound for \(G\)}}

The logarithmic improvement of the lower bound for \(G\) is established under the assumptions
\beeq\label{critical-q}
  q<p,
  \qquad
  pq=2+\frac1q.
\eneq
In particular, \(1<q<2\). Fix
\[
 0<\delta<\frac{1-c}{4}.
\]
Choose \(\chi\in C^\infty_0((0,1))\) such that
\begin{equation*}
 0\leq\chi\leq1,\qquad
 \chi=1\text{ on }[\delta,1-\delta],\qquad
 \supp\chi\subset[\delta/2,1-\delta/2].
\end{equation*}
Choose \(\lambda_0\in(0,1]\) sufficiently small so that
\[
 \frac{\delta}{2\lambda_0}>R.
\]
We may decrease $\lambda_0$ further as needed. For \(0<\lambda\leq\lambda_0\), define
\[
 \Phi_{\lambda,c}(t,r)
 =
 \begin{cases}
 \displaystyle
 \frac{\lambda^{-1/q}}{r}
 \chi\bigl(\lambda(r-ct)\bigr),
 & r>0,\\[2mm]
 0,
 & r=0,
 \end{cases}
\]
and set
\[
 \Phi_c(t,x)
 =
 \int_0^{\lambda_0}
 \Phi_{\lambda,c}(t,|x|)\,d\lambda\ .
\] 
It is easy to verify the following transport identities:
\begin{equation}
\label{eq:Phi-transport-critical-q}
 (\partial_t+c\partial_r)\bigl(r\Phi_{\lambda,c}\bigr)=0,\qquad
 (\partial_t^2-c^2\partial_r^2)\bigl(r\Phi_{\lambda,c}\bigr)=0.
\end{equation}

\begin{lemma}
\label{lem:critical-pairing-q}
Assume \eqref{critical-q}. Then there exist
\(t_1\geq1\) and \(C>0\), independent of \(\vep\), such that
\begin{equation*}
 \int_{\mathbb R^3}u_t(t,x)\Phi_c(t,x)\,dx
 \geq C\vep^{pq}\log(1+t),
 \qquad t_1\leq t<T_\vep.
\end{equation*}
\end{lemma}

\begin{proof}
For $\la\in (0, \la_0)$, define 
\begin{align*}
 H_\lambda(t)
 &=\int_0^\infty U_t(t,r)\bigl(r\Phi_{\lambda,c}\bigr)(t,r)\,dr,
 \\
 F_\lambda(t)
 &=\int_0^\infty(U_t-cU_r)(t,r)\bigl(r\Phi_{\lambda,c}\bigr)(t,r)\,dr.
 \end{align*}
Integration by parts, using \eqref{eq:Phi-transport-critical-q},
gives
\begin{equation*}
 F_\lambda'(t)
 =\int_0^\infty r^{1-p}|V_t(t,r)|^p\bigl(r\Phi_{\lambda,c}\bigr)(t,r)\,dr.
\end{equation*}
Observe that
\[
 2H_\lambda-F_\lambda
 =\int_0^\infty(U_t+cU_r)\bigl(r\Phi_{\lambda,c}\bigr)\,dr.
\]
The characteristic identity
\[
 (\partial_t-c\partial_r)(U_t+cU_r)=r^{1-p}|V_t|^p
\]
implies
\begin{align*}
 (U_t+cU_r)(t,r)
 ={}&(U_t+cU_r)(0,r+ct)\\
 &+\int_0^t
 (r+c(t-s))^{1-p}|V_t(s,r+c(t-s))|^p\,ds.
\end{align*}
On \(\supp\Phi_{\lambda,c}(t,\cdot)\),
\[
 r+ct\geq r-ct\geq\frac{\delta}{2\lambda}
 \geq\frac{\delta}{2\lambda_0}>R\ .
\]
Thus the initial term vanishes, and \(U_t+cU_r\geq0\) there.
The same support condition at \(t=0\) gives
\(F_\lambda(0)=0\). 

These relations will be used to strengthen the lower bound of
\(H_\lambda\), which is the key step toward the improved lower bound for \(G\):
\begin{equation}
\label{eq:core-estimates-critical-q}
 H_\lambda(t)\geq\tfrac12F_\lambda(t)\geq0, \quad F_\lambda'(t)
 =\int_0^\infty r^{1-p}|V_t|^p\bigl(r\Phi_{\lambda,c}\bigr)\,dr,\qquad
 \end{equation}Since \(g_1\geq0\)
and \(g_1\not\equiv0\), there are \(0<a_0<a_1<R\) and \(d_1>0\) such that
\begin{equation}
\label{eq:slow-shell-lower}
 U_t(t,ct+\rho)\geq d_1\vep,\qquad
 a_0\leq\rho\leq a_1,\quad 0\leq t<T_\vep.
\end{equation}
By decreasing \(a_1\), if necessary, we may assume without loss of
generality that
\[
 a_1\leq 2a_0.
\]
We also decrease \(\lambda_0\), if necessary, so that
\[
 \lambda_0^{-1}\geq \frac{2a_0}{\delta}.
\]
Define
\[
 \Omega_\delta=
 \left\{(t,r):(c+\delta)t+a_1\leq r\leq(1-\delta)t,
 \ \ t\geq 2a_0/\delta\right\}.
\]
Fix \((t,r)\in\Omega_\delta\), and set
\[
 s_-=\frac{t-r-a_0}{1+c},\qquad
 s_+=\frac{t-r+a_0}{1-c}.
\]
The definition of \(\Omega_\delta\) implies
\[
 0<s_-<s_+<t.
\]
Moreover, for \(s_-\leq s\leq s_+\),
\[
 |r-t+s|\leq cs+a_0.
\]
On the other hand, since \(s\leq s_+\),
\[
 r+t-s-(cs+a_1)
 \geq r+t-(1+c)s_+-a_1
 >0.
\]
Therefore,
\beeq\label{cap-interval}
 [cs+a_0,cs+a_1]
 \subset
 [|r-t+s|,r+t-s],
 \qquad s_-\leq s\leq s_+.
\eneq

Thus, for these values of \(s\), the whole interval on which the lower
bound for \(U_t\) is known is contained in the spatial integration
interval in the Duhamel formula. The geometry is illustrated in
Figure~\ref{fig:critical-middle-cone}.

\begin{figure}[htbp]
\centering
\resizebox{\textwidth}{!}{%
\begin{tikzpicture}[
  x=0.75cm,
  y=0.77cm,
  >=Stealth,
  line cap=round,
  line join=round,
  every node/.style={font=\large}
]

\pgfmathsetmacro{\Tobs}{10}
\pgfmathsetmacro{\Robs}{7.5}
\pgfmathsetmacro{\cspd}{0.45}
\pgfmathsetmacro{\azero}{0.25}
\pgfmathsetmacro{\aone}{0.55}
\pgfmathsetmacro{\delt}{0.12}
\pgfmathsetmacro{\Jleft}{(\cspd+\delt)*\Tobs+\aone}
\pgfmathsetmacro{\Jright}{(1-\delt)*\Tobs}
\pgfmathsetmacro{\Omegapex}{\aone/(1-2*\delt-\cspd)}

\pgfmathsetmacro{\dref}{\Tobs-\Robs}
\pgfmathsetmacro{\tone}{(\dref-\azero)/(1+\cspd)}
\pgfmathsetmacro{\ttwo}{(\dref+\azero)/(1-\cspd)}

\coordinate (P)      at (\Robs,\Tobs);
\coordinate (Oref)   at (0,\dref);
\coordinate (Dleft)  at (\dref,0);
\coordinate (Dright) at ({\Robs+\Tobs},0);

\coordinate (Pone) at ({\cspd*\tone+\azero},\tone);
\coordinate (Ptwo) at ({\cspd*\ttwo+\azero},\ttwo);

\coordinate (Qone) at ({\cspd*\tone+\aone},\tone);
\coordinate (Qtwo) at ({\cspd*\ttwo+\aone},\ttwo);

\draw[->,very thick]
  (0,0)--({\Robs+\Tobs+1.3},0) node[right] {$\rho$};
\draw[->,very thick]
  (0,0)--(0,{\Tobs+1.3}) node[above] {$s$};

\fill[blue!8]
  (Dleft)--(Dright)--(P)--(Oref)--cycle;

\filldraw[fill=orange!72,fill opacity=0.82,
          draw=orange!90!black,draw opacity=1,very thick]
  ({(1-\delt)*\Omegapex},\Omegapex)
  --(\Jright,\Tobs)
  --(\Jleft,\Tobs)
  --cycle;

\draw[black!65,densely dashed,very thick]
  (0,0)--(\Tobs,\Tobs);
\node[black!70,rotate=45,above,fill=white,fill opacity=0.9,
      text opacity=1,inner sep=1.5pt]
  at (9.65,9.65) {$\rho=s$};

\draw[blue!70!black,very thick]
  (Dleft)--(Dright);

\draw[blue!70!black,very thick]
  (Dleft)--(Oref)--(P);

\draw[blue!70!black,very thick]
  (Dright)--(P);

\node[blue!70!black,below right]
  at ({0.53*\dref},{0.47*\dref})
  {$\rho=t-r-s$};

\node[blue!70!black,above right]
  at ({0.52*\Robs},
      {\dref+0.52*(\Tobs-\dref)})
  {$\rho=s-t+r$};

\node[blue!70!black,above right]
  at ({0.56*(\Robs+\Tobs)+0.44*\Robs},
      {0.44*\Tobs})
  {$\rho=r+t-s$};

\fill (Oref) circle (2.5pt);
\node[left] at (Oref) {$(0,t-r)$};

\fill (P) circle (2.5pt);
\node[above right,xshift=9pt,yshift=7pt] at (P) {$(r,t)$};

\draw[purple!80!black,very thick]
  (\Jleft,\Tobs)--(\Jright,\Tobs);
\draw[purple!80!black,thick]
  (\Jleft,{\Tobs-0.12})--(\Jleft,{\Tobs+0.12})
  (\Jright,{\Tobs-0.12})--(\Jright,{\Tobs+0.12});
\node[purple!80!black,above]
  at ({\Jleft+0.5},\Tobs) {$J_t$};

\draw[purple!60!black,densely dashed,thin]
  (\Jleft,0)--(\Jleft,\Tobs);
\draw[purple!60!black,densely dashed,thin]
  (\Jright,0)--(\Jright,\Tobs);
\node[purple!80!black,below,align=center]
  at (\Jleft,0)
  {$\scriptstyle(c+\delta)t+a_1$};
\node[purple!80!black,below]
  at (\Jright,0)
  {$\scriptstyle(1-\delta)t$};

\draw[green!55!black,very thick]
  plot[
    domain=0:\Tobs,
    samples=2,
    variable=\s
  ]
  ({\cspd*\s+\azero},{\s});

\draw[green!55!black,densely dashed,very thick]
  plot[
    domain=0:\Tobs,
    samples=2,
    variable=\s
  ]
  ({\cspd*\s+\aone},{\s});

\node[green!45!black,left,xshift=-2pt]
  at ({\cspd*7.8+\azero},7.8)
  {$\rho=cs+a_0$};

\node[green!45!black,left,xshift=-5pt,
      fill=white,inner sep=1pt]
  at ({\cspd*9.8+\aone},9.8)
  {$\rho=cs+a_1$};

\fill[green!35,opacity=0.72]
  (Pone)--(Qone)--(Qtwo)--(Ptwo)--cycle;

\draw[green!55!black,thick]
  (Pone)--(Qone)--(Qtwo)--(Ptwo)--cycle;

\node[orange!90!black,fill=white,fill opacity=0.92,text opacity=1,
      inner sep=2pt]
  at (5.65,7.4) {$\Omega_\delta$};

\draw[dashed,gray,thick]
  (0,\tone)--(Pone);

\draw[dashed,gray,thick]
  (0,\ttwo)--(Ptwo);

\fill[red!75!black] (Pone) circle (3.2pt);
\node[red!75!black,below left,align=right,
      fill=white,fill opacity=0.88,text opacity=1,inner sep=1.5pt]
  at (Pone)
  {$\displaystyle s_-=\frac{t-r-a_0}{1+c}$};

\fill[red!75!black] (Ptwo) circle (3.2pt);
\node[red!75!black,above left,align=right,
      fill=white,fill opacity=0.88,text opacity=1,inner sep=1.5pt]
  at (Ptwo)
  {$\displaystyle s_+=\frac{t-r+a_0}{1-c}$};
\end{tikzpicture}
}
\caption{The reflected fast-wave Duhamel region at \((t,r)\) (blue),
the intermediate cone \(\Omega_\delta\) (orange), and the portion of
the slow outgoing shell used in the lower bound (green).
The two intersection times are \(s_-\) and \(s_+\). The horizontal
axis is the radial variable.}
\label{fig:critical-middle-cone}
\end{figure}
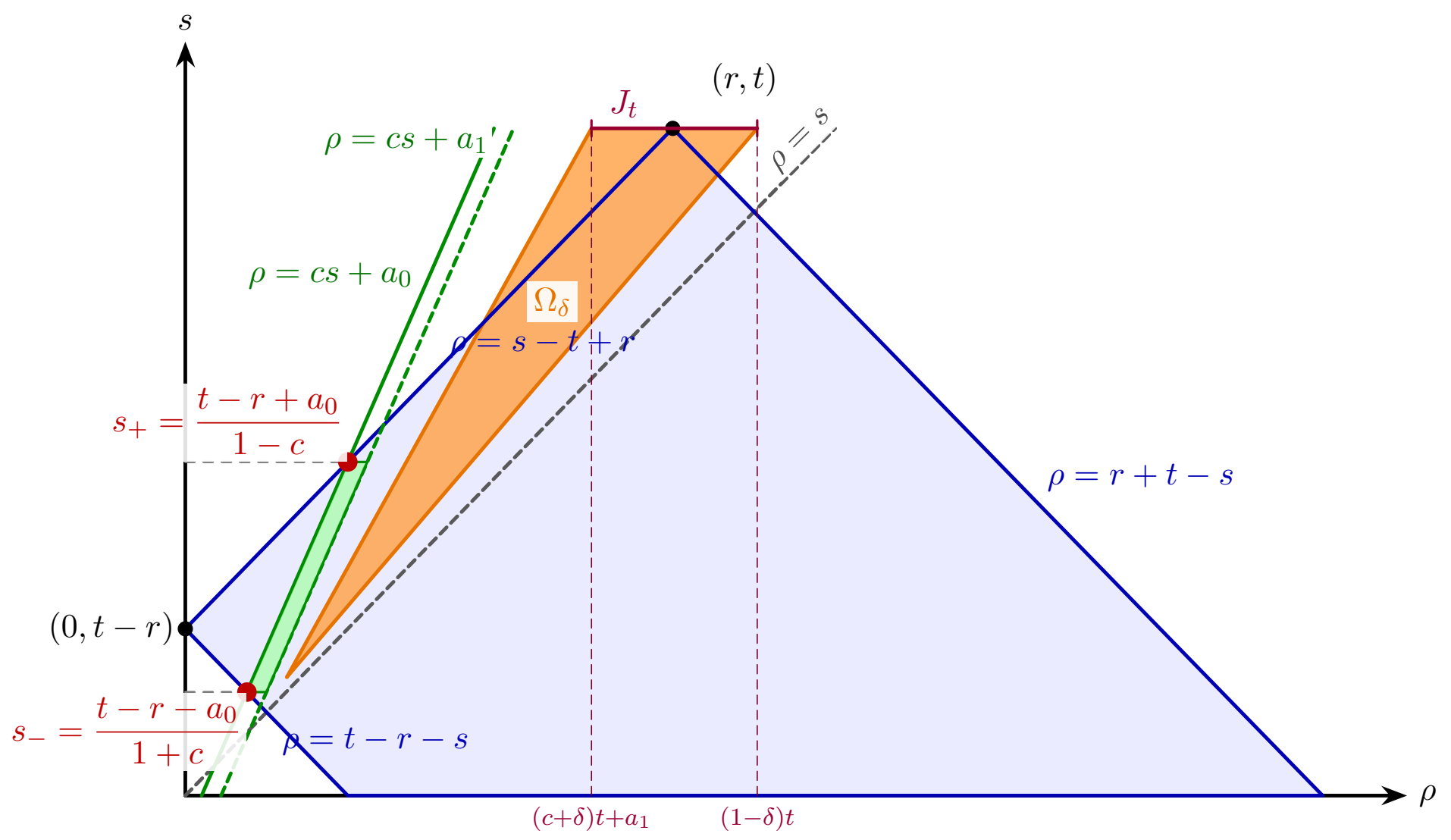
For \((t,r)\in\Omega_\delta\), by the Duhamel formula, \eqref{cap-interval}, and
\eqref{eq:slow-shell-lower}, we have
\begin{align*}
 V(t,r)
 &=
 V^0(t,r)
 +\frac12\int_0^t
   \int_{|r-t+s|}^{r+t-s}
   \rho^{1-q}|U_t(s,\rho)|^q\,d\rho\,ds
\\
 &\geq
 \frac12\int_{s_-}^{s_+}
   \int_{cs+a_0}^{cs+a_1}
   \rho^{1-q}|U_t(s,\rho)|^q\,d\rho\,ds
\\
 &\gtrsim
 \vep^q\int_{s_-}^{s_+}
   \int_{cs+a_0}^{cs+a_1}
   \rho^{1-q}\,d\rho\,ds
\\
 &\gtrsim
 \vep^q t^{1-q}(s_+-s_-)
 \gtrsim
 \vep^q t^{2-q},
\end{align*}
where we have used
\[
 \rho\simeq t, \quad s_+-s_-
 =
 \frac{2c(t-r)+2a_0}{1-c^2}
 \simeq t\ , \quad (t, r)\in \Omega_\delta\ .
\]
 Therefore,
\begin{equation}
\label{eq:V-middle-lower}
 V(t,r)\gtrsim\vep^q t^{2-q},
 \qquad (t,r)\in\Omega_\delta.
\end{equation}
We next convert the pointwise lower bound in
\eqref{eq:V-middle-lower} into a lower bound for the nonlinear source
at the spatial scale \(L\). Set
\[
 L_0=\frac{2a_0}{\delta}.
\]
For \(L_0\leq L<T_\vep\), let
\[
 J_L=[(c+\delta)L+a_1,(1-\delta)L].
\]
By the definition of \(\Omega_\delta\), we have
\[
 \{L\}\times J_L\subset\Omega_\delta,
\]
and \(|J_L|\simeq L\). Hence, by \eqref{eq:V-middle-lower},
\begin{equation*}
 \int_{J_L}V(L,r)\,dr
 \gtrsim\vep^q L^{3-q}.
\end{equation*}
Since \(V(0,r)=0\), the fundamental theorem of calculus gives
\[
 \int_0^L\int_{J_L}|V_t(s,r)|\,dr\,ds
 \geq
 \int_{J_L}V(L,r)\,dr
 \gtrsim\vep^q L^{3-q}.
\]

To recover the weighted source term appearing in
\eqref{eq:core-estimates-critical-q}, we apply H\"older's inequality:
\begin{align*}
 \int_0^L\int_{J_L}|V_t|\,dr\,ds
 &\leq
 \left(
 \int_0^L\int_{J_L}
 r^{1-p}|V_t|^p\,dr\,ds
 \right)^{1/p}
 \left(
 \int_0^L\int_{J_L}r\,dr\,ds
 \right)^{1/p'}.
\end{align*}
Since
\[
 \int_0^L\int_{J_L}r\,dr\,ds\simeq L^3,
\]
we obtain
\begin{equation}
\label{eq:source-mass-critical-q}
 \int_0^L\int_{J_L}
 r^{1-p}|V_t|^p\,dr\,ds
 \gtrsim
 \vep^{pq}L^{3-pq}
 =
 \vep^{pq}L^{1/q'},
\end{equation}
where \(pq=2+1/q\) was used in the last equality.

Then, for \(0<\lambda\leq\lambda_0\), taking \(L=\lambda^{-1}\)
in \eqref{eq:source-mass-critical-q} gives
\begin{equation}
\label{eq:source-mass-lambda-critical-q}
 \int_0^{1/\lambda}\int_{J_{\lambda^{-1}}}
 r^{1-p}|V_t(s,r)|^p\,dr\,ds
 \gtrsim
 \vep^{pq}\lambda^{-1/q'}.
\end{equation}
For \(0\leq s\leq\lambda^{-1}\) and
\(r\in J_{\lambda^{-1}}\),
\[
 \frac{\delta}{\lambda}+a_1
 \leq r-cs
 \leq\frac{1-\delta}{\lambda}.
\]Thus \(\chi(\lambda(r-cs))=1\) throughout this integration region, and
\begin{equation*}
 \bigl(r\Phi_{\lambda,c}\bigr)(s,r)
 =\lambda^{-1/q}.
\end{equation*}
For \(t\geq\lambda^{-1}\), using
\eqref{eq:core-estimates-critical-q} and
\eqref{eq:source-mass-lambda-critical-q}, we find
\begin{align*}
 H_\lambda(t)
 &\geq\tfrac12F_\lambda(t)\\
 &=\tfrac12\int_0^t\int_0^\infty
 r^{1-p}|V_t(s,r)|^p\bigl(r\Phi_{\lambda,c}\bigr)(s,r)\,dr\,ds\\
 &\gtrsim\lambda^{-1/q}
 \int_0^{1/\lambda}\int_{J_{\lambda^{-1}}}
 r^{1-p}|V_t(s,r)|^p\,dr\,ds\\
 &\gtrsim\vep^{pq}\lambda^{-1}.
\end{align*}
In particular,
\begin{equation*}
 H_\lambda(t)\gtrsim\vep^{pq}\lambda^{-1},
 \qquad t\geq\lambda^{-1},\quad0<\lambda\leq\lambda_0.
\end{equation*}
Then, for \(t>t_1=\lambda_0^{-2}+2\), we have 
\begin{align*}
 \int_{\mathbb R^3}u_t(t,x)\Phi_c(t,x)\,dx
 &=4\pi\int_0^{\lambda_0}H_\lambda(t)\,d\lambda\\
 &\geq4\pi\int_{1/t}^{\lambda_0}H_\lambda(t)\,d\lambda\\
 &\gtrsim
 \vep^{pq}\int_{1/t}^{\lambda_0}\frac{d\lambda}{\lambda}
 \ \gtrsim\vep^{pq}\log(1+t).
\end{align*}
This proves the lemma.
\end{proof}

\begin{lemma}[\(L^{q'}\) estimate for \(\Phi_c\)]
\label{lem:Lq-Phi-c}
There exists \(C>0\), independent of \(\vep\), such that
\begin{equation*}
 \int_{|x|\leq t+R}\Phi_c(t,x)^{q'}\,dx
 \leq C(1+t)^{2-q'}\log(2+t),
 \qquad t\geq 1.
\end{equation*}
\end{lemma}

\begin{proof}
Set \(y=r-ct\). On the support of \(\Phi_c\), we have
\[
 y\geq\frac{\delta}{2\lambda_0}>0.
\]
Using the change of variables \(z=\lambda y\), we obtain
\begin{align*}
 \Phi_c(t,r)
 &=\frac{1}{r}\int_0^{\lambda_0}
   \lambda^{-1/q}\chi(\lambda y)\,d\lambda\\
 &=\frac{y^{-1/q'}}{r}
   \int_0^{\lambda_0 y}z^{-1/q}\chi(z)dz\les \frac{y^{-1/q'}}{r}\int_0^{1}z^{-1/q}\chi(z)dz\les \frac{y^{-1/q'}}{r}\ .
\end{align*}
Using polar
coordinates and \(r=ct+y\), we have
\begin{align*}
 \int_{|x|\leq t+R}\Phi_c(t,x)^{q'}\,dx
 &\lesssim
 \int_{\delta/(2\lambda_0)}^{(1-c)t+R}
 (ct+y)^{2-q'}y^{-1}\,dy\\
 &\lesssim
 (1+t)^{2-q'}
 \int_{\delta/(2\lambda_0)}^{(1-c)t+R}\frac{dy}{y}\\
 &\lesssim
 (1+t)^{2-q'}\log(1+t).
\end{align*}
If the integration interval is empty, the desired estimate is
immediate.
\end{proof}

Applying the preceding pairing estimate together with the
\(L^{q'}\)-bound for \(\Phi_c\), and using finite speed of propagation and
H\"older's inequality, we obtain
\[
 \left(\int_{\mathbb R^3}u_t\Phi_c\,dx\right)^q
 \leq G'(t)
 \left(\int_{|x|\leq t+R}\Phi_c^{q'}\,dx\right)^{q-1}.
\]
Hence, by Lemmas~\ref{lem:critical-pairing-q} and
\ref{lem:Lq-Phi-c},
\begin{align*}
 G'(t)
 \gtrsim
 \frac{\vep^{pq^2}(\log(1+t))^q}
      {((1+t)^{2-q'}\log(2+t))^{q-1}}
 \gtrsim
 \vep^{pq^2}(1+t)^{2-q}\log(1+t).
\end{align*}
Integrating over \([t/2,t]\), we obtain
\[
  G(t)
  \gtrsim
  \vep^{pq^2}t^{3-q}\log(1+t)
\]
for all sufficiently large \(t\). This proves
\eqref{eq:critical-log-lower-G} in Lemma \ref{lem:critical-log-improvement}.

\appendix

\section{Scaling-Critical Decay and Its Sharpness}
\label{app:scaling-decay}

We briefly recall the scaling underlying the decay assumptions in
Theorem~\ref{thm:global-critical-decay}. For \(\rho>0\), set
\[
  u_\rho(t,x)=\rho^\alpha u(\rho t,\rho x),
  \qquad
  v_\rho(t,x)=\rho^\beta v(\rho t,\rho x).
\]
Requiring \((u_\rho,v_\rho)\) to satisfy the same system as \((u,v)\)
gives
\[
  \alpha+2=p(\beta+1),
  \qquad
  \beta+2=q(\alpha+1).
\]
Thus
\[
  \alpha_*(p,q)
  =
  \frac{p+2-pq}{pq-1},
  \qquad
  \beta_*(p,q)
  =
  \frac{q+2-pq}{pq-1},
\]
which are precisely the scaling indices appearing in
\eqref{eq:scaling-indices}. Accordingly, the scaling-critical
far-field decay rates for the radial velocity profiles are
\[
  rg_1(r)\simeq(1+r)^{-\alpha_*},
  \qquad
  rg_2(r)\simeq(1+r)^{-\beta_*}.
\]

We now show that both decay rates are sharp. More precisely, for either
component, if its initial velocity profile decays more slowly than the
corresponding scaling-critical rate by any positive power, then finite-time
blow-up occurs for a suitable class of initial data.

\begin{theorem}[Sharpness of the scaling-critical decay]
\label{thm:finite-blowup}
Assume that \(p,q>1\) and \(\de>0\).
Suppose that \(f_1\) and \(f_2\) are compactly supported and that, for some
\(R_0>0\), either
\[
 rg_1(r)\geq(1+r)^{-\alpha_*+\de}, \qquad
g_2\geq0, \quad r\geq R_0\ ,
\]
or
\[
 g_1\geq0, \qquad rg_2(r)\geq(1+r)^{-\beta_*+\de}, \quad r\geq R_0\ .
\]
Then there exist constants \(C_\de>0\) and
\(\vep_\de>0\) such that, for \(0<\vep\leq\vep_\de\), the
corresponding radial solution of \eqref{eq:main-system} cannot exist globally and
\[
  T_\vep\leq C_\de\vep^{-1/\de}.
\]
\end{theorem}

\begin{proof}
Without loss of generality, we may assume that
\[
  c_1=1,
  \qquad
  c_2=a>1.
\]
Let \(R>2R_0\), to be chosen later, be sufficiently large so that
\[
  \supp f_1\cup\supp f_2\subset B_{R/2},
\]
and define
\[
  \Omega
  =
  \left\{
    (t,r):
    t+R\leq r,\quad
    at+\frac{R}{2}\leq r\leq at+R
  \right\}.
\]

We first consider the case
\[
  g_1\geq0,
  \qquad
  rg_2(r)\geq(1+r)^{-\beta_*+\de},
  \quad r\geq R_0.
\]
Since \(f_2\) is compactly supported and
\(r-at\geq R/2\) on \(\Omega\), the terms involving \(f_2\)
in the free solution vanish. Hence, for \((t,r)\in\Omega\),
\begin{align}
\label{eq:appendix-Vt-initial}
  V_t(t,r)\geq &V_t^0(t,r)
  =
  \frac{\vep}{2}
  \left[
    (r+at)g_2(r+at)
    +(r-at)g_2(r-at)
  \right]\\
  \gtrsim&
  \vep R^{-\beta_*+\de}\ .\nonumber
\end{align}
Moreover, since \(r-t\geq R\) on \(\Omega\), the terms involving
\(f_1\) in \(U_t^0\) vanish. Since \(g_1\geq0\), it follows that
\[
  U_t^0(t,r)\geq0,
  \qquad (t,r)\in\Omega.
\]
By \eqref{eq:appendix-Vt-initial} and the Duhamel formula for \(U_t\), we obtain
\begin{align*}
  U_t(t,r)
  &\gtrsim
  \vep^p R^{-p\beta_*+p\de}
  \int_{\frac{r-t-R}{a-1}}^t
    (r-t+s)^{1-p}\,ds
\\
  &\gtrsim
  \vep^p R^{-p\beta_*+p\de}
  r^{1-p}(at-r+R)\ , \qquad (t,r)\in\Omega.
\end{align*}
Hence the iteration may be started with
\[
  U_t(t,r)
  \geq
  A_0
  r^{B_0}
  (r-t-R)^{C_0}
  (at-r+R)^{D_0},
  \qquad
  (t,r)\in\Omega,
\]
where
\[
  A_0\gtrsim
  \vep^p R^{-p\beta_*+p\de},
  \qquad
  B_0=1-p,
  \qquad
  C_0=0,
  \qquad
  D_0=1.
\]
Suppose inductively that
\[
  U_t(t,r)
  \geq
  A_n
  r^{B_n}
  (r-t-R)^{C_n}
  (at-r+R)^{D_n},
  \qquad
  (t,r)\in\Omega.
\]
Using this lower bound in the equation for \(V\), we obtain
\begin{align*}
  V_t(t,r)
  \geq{}&
  A_n^q
  \int_{\frac{at-r+R}{a-1}}^t
    (r-at+as)^{1-q}
    (r-at+as-s-R)^{qC_n}
\\
  &\hspace{2.5cm}\times
    (at-r+R)^{qD_n}
    (r-at+as)^{qB_n}
  \,ds
\\
  \geq{}&
  A_n^q
  r^{qB_n-q+1}
  (at-r+R)^{qD_n}
\\
  &\hspace{1.5cm}\times
  \int_{\frac{at-r+R}{a-1}}^t
    \bigl(
      r-at+(a-1)s-R
    \bigr)^{qC_n}
  \,ds
\\
  ={}&
  \frac{A_n^q}
       {(a-1)(qC_n+1)}
  r^{qB_n-q+1}
  (at-r+R)^{qD_n}
  (r-t-R)^{qC_n+1}.
\end{align*}
Inserting this estimate into the equation for \(U\) yields
\begin{align*}
  U_t(t,r)
  \geq{}&
  \frac{1}{(a-1)(pqD_n+1)}
  \left(
    \frac{A_n^q}
         {(a-1)(qC_n+1)}
  \right)^p
  r^{1-p+p(qB_n-q+1)}
\\
  &\times
  (r-t-R)^{p(qC_n+1)}
  (at-r+R)^{pqD_n+1}.
\end{align*}
Thus we may take
\[
\begin{aligned}
  A_{n+1}
  &=
  \frac{1}{(a-1)(pqD_n+1)}
  \left(
    \frac{A_n^q}
         {(a-1)(qC_n+1)}
  \right)^p,
\\
  B_{n+1}
  &=
  1-p+p(qB_n-q+1),
\\
  C_{n+1}
  &=
  p(qC_n+1),
\\
  D_{n+1}
  &=
  pqD_n+1.
\end{aligned}
\]
Solving these recursions gives
\[
  B_n
  =
  1-p(pq)^n\ ,\quad
  C_n
  =
  \frac{p\bigl((pq)^n-1\bigr)}{pq-1},
  \quad
  D_n
  =
  (pq)^n
  +
  \frac{(pq)^n-1}{pq-1}\ .
\]
Moreover,
\[
  A_n
  =
  A_0^{(pq)^n}
  \prod_{j=0}^{n-1}
  \left[
    (a-1)^{p+1}
    (pqD_j+1)
    (qC_j+1)^p
  \right]^{-(pq)^{n-1-j}}.
\]
Since
\[
  C_j+D_j\lesssim(pq)^j,
\]
the product contributes only an exponential factor of order
\((pq)^n\). Hence there exists \(c_a>0\), independent of \(n\),
\(R\), and \(\vep\), such that
\[
  A_n
  \geq
  \bigl(c_aA_0\bigr)^{(pq)^n}.
\]
Choose
\[
  r_*(t)
  =
  \frac{a+1}{2}t+R.
\]
Then
\[
  r_*(t)-t-R
  =
  at-r_*(t)+R
  =
  \frac{a-1}{2}t,
\]
and
\[
  r_*(t)-at
  =
  R-\frac{a-1}{2}t.
\]
Thus \(r_*(t)\in\Omega\) whenever
\[
  0\leq t\leq\frac{R}{a-1}.
\]
We restrict ourselves to
\[
  I_R
  =
  \left[
    \frac{R}{2(a-1)},
    \frac{R}{a-1}
  \right].
\]
For \(t\in I_R\),
\[
  t\simeq_a R,
  \qquad
  r_*(t)\simeq_a R,
\]
and
\[
  r_*(t)-t-R
  \simeq_a R,
  \qquad
  at-r_*(t)+R
  \simeq_a R.
\]
Substituting \(r=r_*(t)\) into the iterative lower bound gives
\[
\begin{aligned}
  U_t(t,r_*)
  \geq{}&
  A_n
  r_*^{1-p(pq)^n}
  \left(
    \frac{a-1}{2}t
  \right)^{
    \frac{p(q+1)}{pq-1}(pq)^n
    -
    \frac{p+1}{pq-1}
  }.
\end{aligned}
\]
Since
\[
  \beta_*+1
  =
  \frac{q+1}{pq-1},
\]
we rewrite this as
\[
\begin{aligned}
  U_t(t,r_*)
  \geq{}&
  r_*
  \left(
    \frac{a-1}{2}t
  \right)^{-\frac{p+1}{pq-1}}
  \times
  \left[
    c_aA_0
    r_*^{-p}
    \left(
      \frac{a-1}{2}t
    \right)^{p(\beta_*+1)}
  \right]^{(pq)^n}.
\end{aligned}
\]
For \(t\in I_R\),
\[
\begin{aligned}
  A_0
  r_*^{-p}
  \left(
    \frac{a-1}{2}t
  \right)^{p(\beta_*+1)}
  &\gtrsim
  \vep^p
  R^{-p\beta_*+p\de}
  R^{-p}
  R^{p\beta_*+p}
\\
  &=
  \vep^pR^{p\de}.
\end{aligned}
\]
Now choose
\[
  R
  =
  K\vep^{-1/\de},
\]
where \(K>0\) is sufficiently large. For \(\vep>0\) sufficiently small, we also have \(R>2R_0\) and
\[
  \supp f_1\cup\supp f_2\subset B_{R/2},
\]
as required above. Moreover,
\[
  \vep^pR^{p\de}
  =
  K^{p\de}.
\]
Choosing \(K\) sufficiently large, we obtain
\[
  c_aA_0
  r_*^{-p}
  \left(
    \frac{a-1}{2}t
  \right)^{p(\beta_*+1)}
  >1,
  \qquad
  t\in I_R.
\]
Letting \(n\to\infty\) in the iterative lower bound then gives
\[
  U_t(t,r_*)=+\infty,
\]
which contradicts the regularity of the solution. Consequently,
\[
  T_\vep
  \leq
  C_\de R
  \leq
  C_\de\vep^{-1/\de}.
\]
This proves the stated lifespan bound under the first alternative.

We now consider the case
\[
  g_2\geq0,
  \qquad
  rg_1(r)\geq(1+r)^{-\alpha_*+\de},
  \quad r\geq R_0.
\]
We use the same region \(\Omega\) and restrict from the outset to
\[
  0\leq t\leq\frac{R}{a-1}.
\]
Then, for \((t,r)\in\Omega\),
\[
  R\leq r-t\leq2R.
\]
Since \(f_1\) is compactly supported, the terms involving \(f_1\)
in the free solution vanish. Hence
\[
  U_t(t,r)
  \geq
  U_t^0(t,r)
  =
  \frac{\vep}{2}
  \left[
    (r+t)g_1(r+t)
    +(r-t)g_1(r-t)
  \right]
  \gtrsim
  \vep R^{-\alpha_*+\de}.
\]
Moreover, since \(f_2\) is compactly supported and
\(r-at\geq R/2\) on \(\Omega\), the terms involving \(f_2\)
in \(V_t^0\) vanish. Since \(g_2\geq0\), it follows that
\[
  V_t^0(t,r)\geq0,
  \qquad (t,r)\in\Omega.
\]
Thus the same iteration as above may be started with
\[
  U_t(t,r)
  \geq
  A_0
  r^{B_0}
  (r-t-R)^{C_0}
  (at-r+R)^{D_0},
\]
where
\[
  A_0\gtrsim
  \vep R^{-\alpha_*+\de},
  \qquad
  B_0=C_0=D_0=0.
\]
The recurrence relations for \(A_n,B_n,C_n,D_n\) are unchanged. With
these initial values,
\[
  B_n
  =
  1-(pq)^n,
  \qquad
  C_n
  =
  \frac{p\bigl((pq)^n-1\bigr)}{pq-1},
  \qquad
  D_n
  =
  \frac{(pq)^n-1}{pq-1}.
\]
Moreover, as above,
\[
  A_n
  \geq
  \bigl(c_aA_0\bigr)^{(pq)^n}.
\]
Taking again
\[
  r_*(t)
  =
  \frac{a+1}{2}t+R,
  \qquad
  t\in I_R,
\]
we obtain
\[
\begin{aligned}
  U_t(t,r_*)
  \geq{}&
  r_*
  \left(
    \frac{a-1}{2}t
  \right)^{-\frac{p+1}{pq-1}}
  \left[
    c_aA_0
    r_*^{-1}
    \left(
      \frac{a-1}{2}t
    \right)^{\frac{p+1}{pq-1}}
  \right]^{(pq)^n}.
\end{aligned}
\]
Since
\[
  \alpha_*+1
  =
  \frac{p+1}{pq-1},
\]
and
\[
  t\simeq_a R,
  \qquad
  r_*(t)\simeq_a R,
  \qquad t\in I_R,
\]
we have
\[
\begin{aligned}
  A_0
  r_*^{-1}
  \left(
    \frac{a-1}{2}t
  \right)^{\alpha_*+1}
  &\gtrsim
  \vep
  R^{-\alpha_*+\de}
  R^{-1}
  R^{\alpha_*+1}
\\
  &=
  \vep R^\de.
\end{aligned}
\]
Choosing again
\[
  R
  =
  K\vep^{-1/\de},
\]
with \(K>0\) sufficiently large, gives
\[
  \vep R^\de
  =
  K^\de.
\]
Hence, for \(K\) sufficiently large,
\[
  c_aA_0
  r_*^{-1}
  \left(
    \frac{a-1}{2}t
  \right)^{\alpha_*+1}
  >1,
  \qquad
  t\in I_R.
\]
Letting \(n\to\infty\) in the iterative lower bound again gives
\[
  U_t(t,r_*)=+\infty,
\]
which contradicts the regularity of the solution. Therefore
\[
  T_\vep
  \leq
  C_\de\vep^{-1/\de}.
\]
This completes the proof.
\end{proof}

\end{document}